\documentclass[12pt]{amsart}
\usepackage{amssymb}
\usepackage{amscd}

\newcommand{\Q}{\mathbb{Q}}
\newcommand{\Z}{\mathbb{Z}}
\newcommand{\N}{\mathbb{N}}
\newcommand{\C}{\mathbb{C}}

\newcommand{\CalF}{\mathcal{F}}
\newcommand{\CalT}{\mathcal{T}}
\newcommand{\CalO}{\mathcal{O}}
\newcommand{\CalL}{\mathcal{L}}
\newcommand{\CalV}{\mathcal{V}}

\newcommand{\Ind}{\mathrm{Ind}}

\newcommand{\im}{\mathrm{im}\,}

\newcommand{\rmprod}{\mathrm{prod}}
\newcommand{\rmdet}{\mathrm{det}}
\newcommand{\lcm}{\mathrm{lcm}}

\newcommand{\unz}{\underline{\zeta}}
\newcommand{\bbeta}{\overline{\beta}}
\newcommand{\sr}{{(\frac{s}{r})}}

\newcommand{\isomto}{\overset{\sim}{\rightarrow}}
\newcommand{\lng}{\langle}
\newcommand{\rng}{\rangle}

\numberwithin{equation}{section}
\newtheorem{theorem}{Theorem}[section]
\newtheorem{corollary}[theorem]{Corollary}
\newtheorem{lemma}[theorem]{Lemma}
\newtheorem{proposition}[theorem]{Proposition}

\theoremstyle{definition}
\newtheorem{definition}[theorem]{Definition}

\newtheorem{example}[theorem]{Example}

\newcommand{\bdf}{\begin{definition}}
\newcommand{\edf}{\end{definition}\noindent}
\newcommand{\bex}{\begin{example}}
\newcommand{\eex}{\end{example}\noindent}
\newcommand{\bpr}{\begin{proposition}}
\newcommand{\epr}{\end{proposition}}
\newcommand{\blm}{\begin{lemma}}
\newcommand{\elm}{\end{lemma}}
\newcommand{\bth}{\begin{theorem}}
\renewcommand{\eth}{\end{theorem}}
\newcommand{\bpf}{\begin{proof}}
\newcommand{\epf}{\end{proof}\noindent}
\newcommand{\bcr}{\begin{corollary}}
\newcommand{\ecr}{\end{corollary}\noindent}
\newcommand{\beq}{\begin{equation}}
\newcommand{\eeq}{\end{equation}}
\newcommand{\bes}{\begin{equation*}}
\newcommand{\ees}{\end{equation*}}
\newcommand{\ben}{\begin{enumerate}}
\newcommand{\een}{\end{enumerate}}

\begin{document}
\title[Cohomology representations of the symmetric group]
{Bases for certain cohomology representations of the symmetric group}
\author{Anthony Henderson}
\address{School of Mathematics and Statistics,
University of Sydney, NSW 2006, AUSTRALIA}
\email{anthonyh@maths.usyd.edu.au}
\thanks{This work was supported by Australian Research Council grant DP0344185}
\begin{abstract}
We give a combinatorial description (including explicit differential-form 
bases) for the cohomology groups of the space of $n$ distinct nonzero
complex numbers, with coefficients in rank-one local systems which are of
finite monodromy around the coordinate hyperplanes 
and trivial monodromy around all
other hyperplanes. In the case where the local system is
equivariant for the symmetric group, we write the cohomology groups
as direct sums of inductions of one-dimensional characters of subgroups.
This relies on an equivariant description of the Orlik-Solomon algebras
of full monomial reflection groups
(wreath products of the symmetric group with a cyclic group).
The combinatorial models involved are certain representations of
these wreath products which possess bases indexed by labelled trees.
\end{abstract}
\maketitle
\section{Introduction}
Fix a positive integer $n$. This paper concerns the cohomology of the 
complex hyperplane complement
\[ T(1,n):=\{(z_1,z_2,\cdots,z_n)\in \C^{n}\,|\,
z_i\neq 0,\,\forall i,\, z_i\neq z_j,\,\forall i\neq j\}. \]
Let $(\Omega^\bullet(T(1,n)),d)$ be the differential graded algebra
of regular differential forms on $T(1,n)$. The Orlik-Solomon algebra
$A^\bullet(T(1,n))$ is defined to be the 
$\C$-subalgebra of $\Omega^\bullet(T(1,n))$
generated by the following closed $1$-forms corresponding to the hyperplanes:
\[ \omega_i:= \frac{dz_i}{z_i},\,\forall i,\text{ and }
\omega_{i,j}:=\frac{dz_i-dz_j}{z_i-z_j},\,\forall i<j. \]
A famous result of Brieskorn (valid for all hyperplane complements) states
that the inclusion $(A^\bullet(T(1,n)),0)
\hookrightarrow(\Omega^\bullet(T(1,n)),d)$ is a quasi-isomorphism of
cochain complexes, so $A^p(T(1,n))\cong H^p(T(1,n),\C)$ for all $p$.
Moreover, results of Arnold (in this case) and Orlik and Solomon 
(for general hyperplane complements) give a 
combinatorial description 
of $A^\bullet(T(1,n))$, including the following explicit basis.
We say that a rooted forest with vertex set $\{1,\cdots,n\}$ is
\textbf{rectified} if the root of each tree is its largest vertex, and
the path from each other vertex to the root is an increasing sequence.
Let $\CalF(1,n)^\circ$ be the set of such rectified forests equipped
with a partition of the set of roots into two subsets, called
open and closed. To each $F\in\CalF(1,n)^\circ$ associate the
form $\alpha(F)\in A^\bullet(T(1,n))$ which is the wedge product
of all $\omega_{i,j}$ for edges $i$---$j$ of $F$
and all $\omega_i$ for closed roots $i$ of $F$. (For now, we leave
the order of factors in the wedge product unspecified, so
$\alpha(F)$ is defined only up to sign.) Then 
$\{\alpha(F)\,|\,F\in\CalF(1,n)^\circ\}$ is a basis of $A^\bullet(T(1,n))$
(see Theorem \ref{nbcthm} below).

Now suppose we replace the complex coefficients in this cohomology by
a general rank-one local system on $T(1,n)$. Any element
$\omega=\sum_i a_i\omega_i+\sum_{i<j}a_{i,j}\omega_{i,j}$ of
$A^1(T(1,n))$ determines such a local system $\CalL_\omega$, whose local
sections are the solutions of the differential equation $df+f\omega=0$.
We have obvious isomorphisms 
$\CalL_\omega\otimes\CalL_{\omega'}\cong\CalL_{\omega+\omega'}$ for all
$\omega,\omega'\in A^1(T(1,n))$, and $\CalL_\omega\cong\CalL_0=\C$ for
$\omega\in\Z\{\omega_i,\omega_{i,j}\}$. If 
$\ell_\omega:\Omega^\bullet(T(1,n))\to\Omega^\bullet(T(1,n))$ denotes the
map $\alpha\mapsto\omega\wedge\alpha$, then we have a natural isomorphism
\beq
H^p(T(1,n),\CalL_\omega)\cong
\frac{\ker(d+\ell_{\omega}:\Omega^p(T(1,n))\to\Omega^{p+1}(T(1,n)))}
{\im(d+\ell_{\omega}:\Omega^{p-1}(T(1,n))\to\Omega^{p}(T(1,n)))}.
\eeq
For most values of $\omega$, this cohomology too can be described in terms
of the Orlik-Solomon algebra, by a general result of Esnault,
Schechtman, and Viehweg, sharpened by Schechtman, Terao, and Varchenko.
We say that $\omega=\sum_i a_i\omega_i+\sum_{i<j}a_{i,j}\omega_{i,j}$
is \textbf{resonant} if any of the following holds:
\begin{itemize}
\item $\sum_i a_i+\sum_{i<j}a_{i,j}$ is a nonzero integer; or
\item $\sum_{i<j,\,i,j\in I}a_{i,j}$ is a positive integer, for some
subset $I\subseteq\{1,\cdots,n\}$; or
\item $\sum_{i\in I} a_i+\sum_{i<j,\,i,j\in I}a_{i,j}$ is a 
positive integer, for some proper subset $I\subset\{1,\cdots,n\}$.
\end{itemize}
The following result is a special case of \cite[Theorem 9]{stv}
(see also the survey articles \cite{libgoberyuzvinsky} and \cite{yuzvinsky}).
\bth \label{esvthm}
If $\omega$ is not resonant, the inclusion $(A^\bullet(T(1,n)),\ell_\omega)
\hookrightarrow(\Omega^\bullet(T(1,n)),d+\ell_\omega)$ is a quasi-isomorphism
of cochain complexes, so
\[ \frac{\ker(\ell_{\omega}:A^p(T(1,n))\to A^{p+1}(T(1,n)))}
{\im(\ell_{\omega}:A^{p-1}(T(1,n))\to A^{p}(T(1,n)))}\cong 
H^p(T(1,n),\CalL_\omega). \]
\eth
\noindent
(Note that Brieskorn's result is the special case $\omega=0$.)
Even if $\omega$ itself is resonant, it is often possible to replace it with
a non-resonant element of $\omega+\Z\{\omega_i,\omega_{i,j}\}$, giving
an isomorphic local system. However, there exist $\omega$ such that
all elements of $\omega+\Z\{\omega_i,\omega_{i,j}\}$ are resonant:
for instance, it is easy to see that $\frac{s}{n}\sum_i\omega_i$
has this property when $s\in\Z$ and $1<\gcd(s,n)<n$.
The corresponding local systems cannot be handled by Theorem \ref{esvthm}.

One result of this paper is an explicit basis of
$H^\bullet(T(1,n),\CalL_{\sum a_i\omega_i})$ where all $a_i\in\Q$
(note that the coefficients of all $\omega_{i,j}$'s here are zero).
This basis is given by differential forms attached to the subset
of $\CalF(1,n)^\circ$ consisting of forests which satisfy
\beq \label{integraleqn}
\sum_{i\in T} a_i\in\Z,\text{ for all trees $T$.}
\eeq
(Note that if $\sum_{i=1}^n a_i\not\in\Z$, this subset is empty.)
The associated differential form $\bbeta(F)\in\Omega^\bullet(T(1,n))$ 
is usually not in $A^\bullet(T(1,n))$, because there are two
changes from the above definition of $\alpha(F)$. Firstly, any edge $i$---$j$
which is \textbf{unbreakable}, in the sense that deleting it 
results in a forest
which no longer satisfies \eqref{integraleqn}, contributes to the wedge
product not $\omega_{i,j}$ but rather $\frac{z_j dz_i-z_i dz_j}
{z_i z_j(z_i-z_j)}$. Secondly, the wedge product must be multiplied by
a certain monomial in the $z_i$'s, to be defined in Section 4.
The result (see Theorem \ref{finalthm} and Corollary \ref{finalcor}) 
is that each such $\bbeta(F)$
lies in $\ker(d+\ell_{\sum a_i\omega_i})$, and their images in
$H^\bullet(T(1,n),\CalL_{\sum a_i\omega_i})$ form a basis.
A consequence (Corollary \ref{modulecor}) is that the 
$H^\bullet(T(1,n),\C)$-module $H^\bullet(T(1,n),\CalL_{\sum a_i\omega_i})$
is generated by those $\bbeta(F)$'s which correspond to forests in which
all roots are open and all edges are unbreakable.
See Example \ref{mysteryex} for the special case where all $a_i=\frac{s}{n}$
and $\gcd(s,n)=1$, which can also be handled by a result of Kawahara.

The special property of the local systems $\CalL_{\sum a_i\omega_i}$,
$a_i\in\Q$, which makes such results possible (indeed, easy) is that
they become trivial on pulling back to a covering space of $T(1,n)$
which is \emph{also a hyperplane complement}. Namely, let $r$ be a
positive integer such that $ra_i\in\Z$ for all $i$, and define
\[  T(r,n):= \{(z_1,z_2,\cdots,z_n)\in \C^{n}\,|\,
z_i\neq 0,\,\forall i,\, z_i^r\neq z_j^r,\,\forall i\neq j\}. \]
The map $\varphi:T(r,n)\to T(1,n):(z_1,\cdots,z_n)\mapsto(z_1^r,\cdots,z_n^r)$
is an unramified Galois covering with group $\mu_r^n$,
where $\mu_r$ is the cyclic group consisting of all complex $r$th roots of $1$.
The local system $\CalL_{\sum a_i\omega_i}$ is the direct summand
of the push-forward $\varphi_*(\C)$ corresponding to the character
$(\zeta_1,\cdots,\zeta_n)\mapsto\zeta_1^{-ra_1}\cdots\zeta_n^{-ra_n}$ 
of $\mu_r^n$, so $H^p(T(1,n),\CalL_{\sum a_i\omega_i})$ is isomorphic via
the pull-back $\varphi^*$ to the corresponding isotypic component 
of $H^p(T(r,n),\C)$ (for more details
on this isomorphism, see Proposition \ref{cochainprop}).
The Orlik-Solomon results on $A^\bullet(T(r,n))$ lead easily to
a basis $\{\beta(F)\}$
of this isotypic component, and a short calculation shows that this is
the pull-back of $\{\bbeta(F)\}$. 

In fact I knew
the basis $\{\beta(F)\}$ on $T(r,n)$ first, and was surprised to find
a relatively nice formula for the corresponding basis on $T(1,n)$.
However, one of the referees has suggested an alternative approach which
avoids (explicitly) passing to the covering space, and thus removes
some of the surprise. Namely, since the local system $\CalL_{\sum a_i\omega_i}$
has trivial monodromy around the hyperplanes $z_i=z_j$, one could successively
delete these using the deletion-restriction result of Cohen 
(\cite[Theorem 4]{cohen}), and relate 
$H^\bullet(T(1,n),\CalL_{\sum a_i\omega_i})$ to
$H^\bullet(T(1,m),\CalL_{\sum a_i\omega_{f(i)}})$ for various
$m<n$ and maps $f:\{1,\cdots,n\}\to\{1,\cdots,m\}$; the rectified forests
and the integrality condition \eqref{integraleqn} would reappear.
Essentially, what this means is that one can mimic on the base space
the steps of the proof of the Orlik-Solomon results on the covering space.

If this basis result had been the only goal in view, this paper would have
been much shorter. But after the paper \cite{lehrersolomon} of
Lehrer and Solomon, one has a right to expect of such Orlik-Solomon-style
descriptions that they take into account the action of the relevant symmetry
group. In the present case,
the symmetric group $S_n$ acts on $T(1,n)$ by permuting coordinates, and
each cohomology group $H^p(T(1,n),\CalL_{\sum a_i\omega_i})$ is thus
a representation for the subgroup of $S_n$ which fixes $\sum a_i\omega_i$,
namely 
\[ Z_{(a_i)}:=\{w\in S_n\,|\, a_{w(i)}=a_i,\,\forall i\}. \] 
On the model of \cite[(4.5)]{lehrersolomon}, one should aim
to write this representation as a direct sum of inductions
of one-dimensional characters of subgroups. Corollary \ref{othermainindcor}
below accomplishes this in the case 
where all $a_i$ are equal, so $Z_{(a_i)}=S_n$; having understood
this case, the interested reader will have no trouble imagining the general 
result. The reason for concentrating on this case (apart from 
simplifying the notation) is that it 
arose in \cite{mytorus}, and the motivation for all
this work was to give a better explanation for the following isomorphism
of representations of $S_n$,
which I originally proved by computing characters:
\beq \label{intromysteryeqn}
H^p(T(1,n),\CalL_{\frac{s}{r}\sum\omega_i})\cong
\varepsilon_n\otimes\Ind_{W(r,n/r)}^{S_n}(\rmdet_{n/r}\otimes
H^{p-n+n/r}(T(r,n/r),\C)).
\eeq
Here $\gcd(s,r)=1$ (so the fraction $\frac{s}{r}$ is in lowest terms), 
$\varepsilon_n$ is the sign character of $S_n$,
$W(r,m)$ is the wreath product $\mu_r\wr S_m$, and $\rmdet_{m}$
is the determinant of the natural representation of $W(r,m)$ on $\C^m$.
For the new explanation, see Corollary \ref{mysterycor}.

In order to study these representations of $S_n$, we need a sufficiently
equivariant description of $A^\bullet(T(r,n))$, one which takes into account
the action of $W(r,n)$. Section 3 is
devoted to such a description; in Corollary \ref{mainindcor} we 
write $A^\bullet(T(r,n))$ as a direct sum
of inductions of one-dimensional characters of subgroups of $W(r,n)$,
generalizing a result of Douglass (\cite[(1.1)]{douglass})
in the $r=2$ case. Since the Orlik-Solomon basis is not stable
under $W(r,n)$, our approach is to 
consider the collection of all $W(r,n)$-translates
of basis elements and the linear dependence relations they satisfy.
Being suitably careful with the signs, we will find that these relations
can be described in a uniform combinatorial way, using rooted forests
which are now not necessarily rectified.
This motivates the preparatory Section 2, where we study 
a representation of $W(r,n)$
defined abstractly 
by certain relations between labelled trees, inspired by
\cite{lehrersolomon} and \cite{barcelo}.

\noindent
\textit{Acknowledgements. }I am very grateful to H.~Barcelo, G.~Lehrer,
and the referees
for their helpful comments on this work.
\section{Representations of wreath products on trees}
Fix positive integers $r$ and $n$. As in the introduction, let 
$\mu_r$ denote the group of $r$th complex roots of $1$,
$S_n$ the group of permutations of $\{1,\cdots,n\}$, and $W(r,n)$
the wreath product $\mu_r\wr S_n$, of cardinality
$r^n n!$. In other words, $W(r,n)$ is the semidirect product 
$\mu_r^n \rtimes S_n$, where $S_n$ acts on
$\mu_r^n$ by permuting the factors; we will usually write
its elements in the form $(\zeta_1,\cdots,\zeta_n)w$, where
$\zeta_i\in\mu_r$ and $w\in S_n$. Often we will abbreviate
$(\zeta_1,\cdots,\zeta_n)$ as $\unz$.
We need some one-dimensional characters of $W(r,n)$:
\begin{itemize}
\item $\varepsilon_n$ defined by $\varepsilon_n(\unz w)=\varepsilon_n(w)$,
the extension of the sign character of $S_n$;
\item $\rmprod_n$ defined by
$\rmprod_n(\unz w)=\zeta_1\cdots\zeta_n$;
and 
\item $\rmdet_n:=\varepsilon_n\rmprod_n$.
\end{itemize} 
The last is so named because it is
the determinant arising from the standard representation of $W(r,n)$
on $\C^n$, in which $w$ is represented by the usual permutation matrix and
$\unz$ by the diagonal matrix with entries $(\zeta_1,\cdots,\zeta_n)$.
The image of $W(r,n)$ under this representation is the monomial reflection
group $G(r,1,n)$.

In this section
we study a representation of $W(r,n)$ which has a basis indexed
by certain labelled trees. Here are our conventions.
\bdf
A \textbf{directed tree} on a nonempty finite set $I$ is a directed graph $T$
with vertex set $I$, containing no loops or cycles, 
such that there is exactly one vertex, the \textbf{root},
with out-degree $0$, and every other vertex has out-degree $1$.
\edf
It is clear that these conditions force the underlying graph of $T$
to be a tree, and that all edges are directed `towards the root'.
So by a well known result, the total number of directed trees 
on $I$ is $|I|^{|I|-1}$.
\bdf
We define a \textbf{$\mu_r$-labelled directed tree}
on $I$ to be a directed tree on $I$ equipped with a labelling
of each edge by an element of $\mu_r$. (If $r=1$, the labels are all $1$
and may be neglected.)
If $T$ is a $\mu_r$-labelled directed tree,
we write $i\overset{\eta}{\rightarrow}j$ (with a subscript $T$ when necessary)
to mean that $i,j\in I$ and 
there is an edge from $i$ to $j$ in $T$ whose label is $\eta\in\mu_r$.
Let $\CalT(r,I)$ be the set of $\mu_r$-labelled directed
trees on $I$, and write $\CalT(r,n)$ for $\CalT(r,\{1,\cdots,n\})$.
\edf
Clearly the cardinality of $\CalT(r,n)$ is $(rn)^{n-1}$.
We have a natural action of $W(r,n)$ on $\CalT(r,n)$, in which
$S_n$ acts by permuting the vertex set and $\unz\in\mu_r^n$
acts by multiplying the label of an edge from $i$ to $j$ by
$\zeta_i\zeta_j^{-1}$. More formally, the following equivalences hold:
\beq
\begin{array}{ccccc}
i\overset{\eta}{\longrightarrow}j & \Longleftrightarrow &
i\overset{\ \eta\zeta_i\zeta_j^{-1}}{\longrightarrow} j &
\Longleftrightarrow &
w(i)\overset{\eta}{\longrightarrow} w(j)\\
T&&\unz.T&&w.T
\end{array}
\eeq
for all $\unz\in\mu_r^n$ and $w\in S_n$.
Linearizing, we obtain a representation of $W(r,n)$ on the
vector space $\C\CalT(r,n)$ with $\CalT(r,n)$ as basis.

Now we want to define a quotient $\CalV{(r,n)}$ of $\C\CalT(r,n)$
by imposing some linear relations which the images $[T]$ of
$T\in\CalT(r,n)$ must satisfy. (The motivation for these relations will
become clear when 
we consider differential forms in the next section.) The relations
are of the following two kinds:
\beq \label{firsttreereleqn}
[T_1]+[T_2]=[T_3],
\eeq
if $T_1,T_2,T_3\in\CalT(r,n)$
are identical except for the edges between three vertices $i$, $j$, and $k$,
where we have

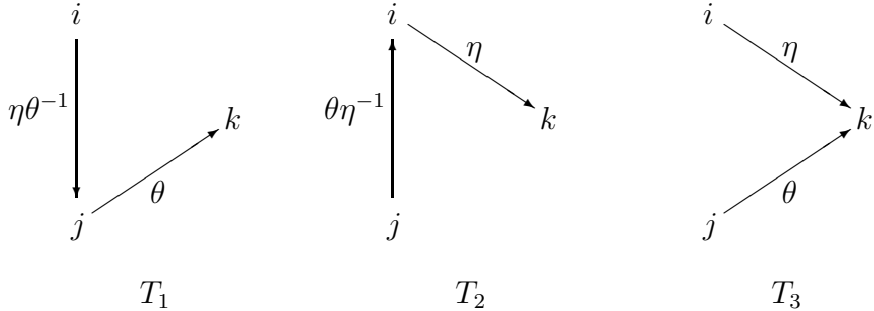
\setlength{\unitlength}{0.7mm}
\begin{picture}(160,70)(0,0)
\put(9,58){$i$}
\put(9,18){$j$}
\put(38,38){$k$}
\put(10,55){\vector(0,-1){30}}
\put(13,22){\vector(3,2){24}}
\put(-3,40){$\eta\theta^{-1}$}
\put(24,24){$\theta$}
\put(22,5){$T_1$}
\put(69,58){$i$}
\put(69,18){$j$}
\put(98,38){$k$}
\put(70,25){\vector(0,1){30}}
\put(73,58){\vector(3,-2){24}}
\put(57,40){$\theta\eta^{-1}$}
\put(84,52){$\eta$}
\put(82,5){$T_2$}
\put(129,58){$i$}
\put(129,18){$j$}
\put(158,38){$k$}
\put(133,58){\vector(3,-2){24}}
\put(133,22){\vector(3,2){24}}
\put(144,24){$\theta$}
\put(144,52){$\eta$}
\put(142,5){$T_3$}
\end{picture}

\noindent
for some $\eta,\theta\in\mu_r$; and
\beq \label{secondtreereleqn}
[T]+[T']=0,
\eeq
if $T,T'\in\CalT(r,n)$
are identical except for the edge between two vertices $i$ and $j$,
where we have:

\setlength{\unitlength}{0.7mm}
\begin{picture}(160,32)(0,23)
\put(27,40){\vector(1,0){26}}
\put(123,40){\vector(-1,0){26}}
\put(22,38){$i$}
\put(56,38){$j$}
\put(92,38){$i$}
\put(126,38){$j$}
\put(38,43){$\eta$}
\put(107,43){$\eta^{-1}$}
\put(38,28){$T$}
\put(107,28){$T'$}
\end{picture}

\noindent
for some $\eta\in\mu_r$. (Note that this implies that
$j$ is the root of $T$ and $i$ is the root of $T'$.)
It is clear that $W(r,n)$ permutes these relations, so we obtain
a representation of $W(r,n)$ on $\CalV{(r,n)}$.

One feature of these relations is that they allow us to express
any tree as a linear combination of those where the edges are always
directed towards the greater vertex, a process of `rectification' which
is implicit in \cite{lehrersolomon} and \cite{barcelo}.
This is encapsulated in the next definition and result.
\bdf
A tree $T\in\CalT(r,n)$ is \textbf{rectified} if 
$i\overset{\eta}{\rightarrow}j$ implies $i<j$. Let
$\CalT(r,n)^\circ$ be the set of rectified trees.
\edf
Note that $\CalT(r,n)^\circ$ has exactly $r^{n-1}(n-1)!$ elements:
the root must be $n$, and for $i<n$, there are $n-i$ possibilities
for the end of the edge which starts at $i$, with $r$ possible labels
for the edge.
\blm \label{spanlemma}
$\CalV{(r,n)}$ is spanned by $\{[T]\,|\,T\in\CalT(r,n)^\circ\}$.
\elm
\bpf
By using relations of type \eqref{secondtreereleqn} repeatedly to shift
the root,
we can express any $[T]$ for $T\in\CalT(r,n)$ as $\pm[T']$ where
$T'$ has root $n$. Now suppose that 
$T\in\CalT(r,n)$ has root $n$ and is unrectified. 
let $a(T)\leq n-1$ be the maximal $i$ such that $T$ has an edge
$i\overset{\eta}{\rightarrow}j$ where $i>j$, and let
$b(T)$ be the length of the path from $a(T)$ to $n$ in $T$.
It suffices to show that
$[T]$ can be written as a linear combination of those $[T']$ where
$T'$ has root $n$ and is either rectified or has a pair
$(a(T'),b(T'))$ which precedes $(a(T),b(T))$ in lexicographic order. 
For this, define $T_2,T_3\in\CalT(r,n)$ so that $T,T_2,T_3$ form a triple
as in \eqref{firsttreereleqn}, where $i=a(T)$, $j$ is the end of the
edge starting at $i$ in $T$, and $k$ is the end of the edge starting at
$j$ in $T$ ($j$ cannot be the root $n$, since by assumption it is less
than $i$). Since $T_2$ and $T_3$ clearly have the same root as $T$,
it suffices to show that they are either
rectified or have an $(a,b)$ pair preceding that of $T$. Now if
$k>a(T)=i$, then the edge starting with $i$ in both $T_2$ and $T_3$
is directed towards the greater vertex, so they are either rectified
or have an $a$-value strictly less than $a(T)$. On the other hand,
if $k<a(T)$, then $a(T_2)=a(T_3)=a(T)$, but $b(T_2)=b(T_3)=b(T)-1$, since
$k$ is a step closer to the root than $j$ is in $T$. 
So in either case we are done.
\epf

Now the representation $\CalV{(1,n)}$ of $S_n$ has a convenient
`concrete' realization. For $T\in\CalT(1,n)$, define the
polynomial $p_T\in\C[z_1,\cdots,z_n]$ by 
$p_T:=\prod_{\substack{i\rightarrow j\\T}}(z_i-z_j)$.
Let $S_n$ act on
$\C[z_1,\cdots,z_n]$ and its quotient field $\C(z_1,\cdots,z_n)$ in
the obvious way, so that $w.z_i=z_{w(i)}$.
\bpr \label{r1prop}
\ben
\item $\exists$ a
linear map $\varrho:\CalV{(1,n)}\to\C(z_1,\cdots,z_n)$
such that $\varrho([T])=\frac{1}{p_T},\ \forall T\in\CalT(1,n)$.
\item The map $\varrho$ is injective and $S_n$-equivariant.
\item $\CalV{(1,n)}$ has basis 
$\{[T]\,|\,T\in\CalT(1,n)^\circ\}$.
\een
\epr
\bpf
To prove (1), it suffices to
show that $\frac{1}{p_{T_1}}+\frac{1}{p_{T_2}}=\frac{1}{p_{T_3}}$
for $T_1,T_2,T_3$ a triple as in \eqref{firsttreereleqn}, and
$\frac{1}{p_T}+\frac{1}{p_{T'}}=0$ for $T,T'$ as in \eqref{secondtreereleqn};
these are both trivial. The equivariance part of (2) is also obvious.
Now by Lemma \ref{spanlemma}, the set in (3)
is a spanning set, so (2) and (3) will both follow from the
statement that $\{\frac{1}{p_T}\,|\,T\in\CalT(1,n)^\circ\}$ is a linearly
independent subset of $\C(z_1,\cdots,z_n)$. We could deduce this from 
Orlik-Solomon theory using \eqref{nbceqn} below, but there is also
a pleasant direct proof by induction on $n$.
The claim is vacuously true when $n=1$, so suppose that
$n\geq 2$ and that 
\beq  \label{linindeqn}
\sum_{T\in\CalT(1,n)^\circ}\frac{c_T}{p_T}=0,
\eeq
for some $c_T\in\C$. Any 
$T\in\CalT(1,n)^\circ$ contains the edge $n-1\to n$, and hence
determines a partition $\{1,\cdots,n\}=A_T\sqcup B_T$, where
$A_T$ consists of all the vertices of $T$ closer to $n-1$ than to $n$,
and $B_T$ consists of all the vertices closer to $n$. 
Fix a partition $\{1,\cdots,n\}=A\sqcup B$,
where $n-1\in A$ and $n\in B$, and define
\[ \CalT(1,n)_{A,B}^\circ=\{T\in\CalT(1,n)^\circ\,|\,A_T=A,B_T=B\}. \]
It suffices to show that $c_T=0$ for all $T\in\CalT(1,n)_{A,B}^\circ$.
Now we can rewrite \eqref{linindeqn} as
\beq 
\sum_{T\in\CalT(1,n)_{A,B}^\circ}
c_T\prod_{\substack{i\to j\\T\\i\leq n-2}}(z_i-z_j)^{-1}
+\sum_{\substack{T\in\CalT(1,n)^\circ\\T\not\in\CalT(1,n)_{A,B}^\circ}}
c_T\prod_{\substack{i\to j\\T\\i\leq n-2}}(z_i-z_j)^{-1}=0,
\eeq
where we have multiplied through by the common factor $z_{n-1}-z_n$.
The terms in the first sum contain only 
factors $(z_i-z_j)^{-1}$ where
$i,j$ are both in $A$ or both in $B$, whereas each of the terms in the
second sum contains at least one factor $(z_i-z_j)^{-1}$ where one of $i,j$ is
in $A$ and the other is in $B$. Therefore the first sum is itself zero:
one way to see this is to make a substitution of variables
$z_a\mapsto z_a+s$ for all $a\in A$ and $z_b\mapsto z_b+t$ for all $b\in B$,
and then send $|s-t|\to\infty$, killing all terms in the second sum while
leaving the first unchanged. We are now reduced to showing that the set
\beq \label{seteqn}
\{\prod_{\substack{i\to j\\T\\i\leq n-2}}(z_i-z_j)^{-1}\,|\,
T\in\CalT(1,n)_{A,B}^\circ\}
\eeq
is linearly independent. But a tree $T\in\CalT(1,n)_{A,B}^\circ$
is constructed from uniquely defined smaller trees $T|_A\in\CalT(1,A)^\circ$
and $T|_B\in\CalT(1,B)^\circ$, by adding the edge from $n-1$ to $n$:
so the set \eqref{seteqn} is exactly 
\[ \{\frac{1}{p_{T'}p_{T''}}\,|\,
T'\in\CalT(1,A)^\circ,T''\in\CalT(1,B)^\circ\}. \] 
By the induction
hypothesis, $\{\frac{1}{p_{T'}}\,|\,T'\in\CalT(1,A)^\circ\}$ is a
linearly independent set in $\C(z_a\,|\,a\in A)$, and
$\{\frac{1}{p_{T''}}\,|\,T''\in\CalT(1,B)^\circ\}$ is a
linearly independent set in $\C(z_b\,|\,b\in B)$, so the claim follows.
\epf

We can now deduce the analogous basis for $\CalV{(r,n)}$ for general $r$.
\bpr \label{indprop}
Define
$Z:=\{\unz=(\zeta_1,\cdots,\zeta_{n})\in\mu_r^n\,|\,\zeta_n=1\}$.
\ben
\item $\CalV(1,n)$ can be identified with the subspace of $\CalV(r,n)$
spanned by $[T]$ for $T\in\CalT(1,n)$.
\item $\CalV{(r,n)}=\bigoplus_{\unz\in Z}\ \unz.\CalV{(1,n)}$.
\item $\CalV{(r,n)}$ has basis $\{[T]\,|\,T\in\CalT(r,n)^\circ\}$.
\item $\dim\CalV{(r,n)}=r^{n-1}(n-1)!$.
\item Embed $\mu_r\times S_n$ in $W(r,n)$ by the map
$(\zeta,w)\mapsto(\zeta,\cdots,\zeta)w$. As a representation of $W(r,n)$,
\bes
\CalV{(r,n)}\cong\Ind_{\mu_r\times S_n}^{W(r,n)}(\CalV{(1,n)}),
\ees
where $\mu_r$ acts trivially on $\CalV{(1,n)}$.
\een
\epr
\bpf
Identifying $\CalT(1,n)$ with the subset
of $\CalT(r,n)$ where all the edge-labels are $1$, we have that for
$\unz\in Z$, 
\beq 
\unz.\CalT(1,n)=\{T\in\CalT(r,n)\,|\,i\overset{\eta}{\rightarrow} j
\Longrightarrow \eta=\zeta_i\zeta_j^{-1}\}.
\eeq
It is clear that for any tree $T\in\CalT(r,n)$,
the edge labels satisfy this rule for a unique $\unz\in Z$,
so $\CalT(r,n)$ is the disjoint
union of these subsets $\unz.\CalT(1,n)$. Moreover,
the relations involve only trees in the same subset, 
which proves (1) and (2). Parts (3) and (4) follow from (2) and Proposition 
\ref{r1prop}. Part (5) follows from (2) and the fact that $Z$ is a set of coset
representatives for $W(r,n)/(\mu_r\times S_n)$.
\epf

There is another easy-to-describe basis of $\CalV{(r,n)}$.
Let $T_0\in\CalT(1,n)$ be the tree (a chain)
in which $i\to i+1$ for all $i\leq n-1$.
For $T\in\CalT(r,n)$, $\unz\in Z$ and $w\in S_n$, we say that $(\unz,w)$ 
\textbf{refines} $T$
if $i\overset{\eta}{\to} j$ in $T$ 
implies $w^{-1}(i)<w^{-1}(j)$ and $\eta=\zeta_i\zeta_j^{-1}$; 
in other words, $T\in\unz.\CalT(1,n)$ and the
total order on $\{1,\cdots,n\}$ determined by the chain $w.T_0$ refines
the partial order determined by $T$.
\bpr \label{unbranchingprop}
For all $T\in\CalT(r,n)$, the equation
\[ [T]=\sum_{\substack{(\unz,w)\in Z\times S_n\\(\unz,w)\textrm{ refines }T}} 
[\unz w.T_0] \]
holds in $\CalV{(r,n)}$.
\epr
\bpf
Let $f(T)$ denote the right-hand side. If $T_1,T_2,T_3\in\CalT(r,n)$ are as in
\eqref{firsttreereleqn}, it is clear that the set of $(\unz,w)$ which refine
$T_3$ is the disjoint union of the corresponding sets for $T_1$ and
$T_2$, and hence $f(T_3)=f(T_1)+f(T_2)$. Now certainly $[T]=f(T)$
if $T$ is a chain, i.e.\ $T=\unz w.T_0$ for some $(\unz,w)$. If $T$
is not of this form, then it forms part of a triple $(T_1,T_2,T)$
as in \eqref{firsttreereleqn}, where $T_1$ and $T_2$ strictly precede
$T$ in the lexicographic ordering according to the $n$-tuple
$(v_0,\cdots,v_{n-1})$, where $v_d$ is the number of vertices at distance
$d$ from the root. The result follows by induction.
\epf
\bcr \label{newbasiscor}
$\CalV{(r,n)}$ has basis $\{[\unz w.T_0]\,|\,\unz\in Z, w\in S_{n-1}\}$.
\ecr
\bpf
If $T\in\CalT(r,n)^\circ$,
the partial order corresponding to $T$ has unique maximal element $n$,
so if $(\unz,w)$ refines $T$, $w$ must lie in $S_{n-1}$.
So (3) of Proposition \ref{indprop} and Proposition \ref{unbranchingprop}
imply that the stated set spans $\CalV{(r,n)}$; it must be a basis
by (4) of Proposition \ref{indprop}.
\epf
These results offer an alternative way of defining $\CalV{(r,n)}$.
Consider the quotient $\widetilde{\CalV}{(r,n)}$ of $\C\CalT(r,n)$
obtained by imposing only the relations \eqref{firsttreereleqn}, and
not \eqref{secondtreereleqn}. The above arguments show
that the images of the chains in $\CalT(r,n)$
form a basis of $\widetilde{\CalV}{(r,n)}$. We can then regard
$\CalV{(r,n)}$ as the quotient of the formal span of these chains
by the relations \eqref{secondtreereleqn}, where $[T]$ and $[T']$ are
interpreted by means of Proposition
\ref{unbranchingprop}. Another way of saying this is that we can
write explicit generators for the left ideal $I$ of $\C W(r,n)$
defined by $I=\{x\in\C W(r,n)\,|\,x.[T_0]=0\}$: namely,
\[ (\zeta,\cdots,\zeta)-1,\text{ for }\zeta\in\mu_r,\text{ and}\negthickspace
\sum_{w\in (S_d\times S_{n-d})\setminus S_n}\negthickspace\negthickspace 
w,\text{ for }1\leq d\leq n-1. \]
(Here $(S_d\times S_{n-d})\setminus S_n$ denotes the set of minimal-length
left coset representatives for the Young subgroup $S_d\times S_{n-d}$.)
Then $x\mapsto x.[T_0]$ gives an isomorphism $\C W(r,n)/I\isomto\CalV{(r,n)}$
of representations of $W(r,n)$. This point of view will not be used in the
remainder of the paper.

Lehrer and Solomon have effectively already studied the representation
$\CalV{(1,n)}$. Indeed, in \cite[Section 3]{lehrersolomon} they consider
the top-degree component $A^{n-1}(M(n))$ of the Orlik-Solomon algebra
of the hyperplane complement
\[ M(n):=\{(z_1,\cdots,z_n)\in\C^n\,|\,z_i\neq z_j,\,\forall i\neq j\}. \]
This component has a `no-broken-circuit'
basis $\{\alpha_T\,|\,T\in\CalT(1,n)^\circ\}$
where $\alpha_T:=\bigwedge_{i\to j}\frac{dz_i-dz_j}{z_i-z_j}$, with
the factors of the wedge product ordered
by increasing $i$. But by an elementary calculation,
\beq \label{nbceqn}
\alpha_T=\frac{1}{p_T}\alpha,\text{ where }
\alpha:=\sum_{i=1}^n (-1)^{n-i}dz_1\wedge\cdots\wedge\widehat{dz_i}\wedge
\cdots\wedge dz_n.
\eeq
(The hat denotes omission of that factor.) 
So we have a linear isomorphism $\CalV{(1,n)}\isomto
A^{n-1}(M(n)):x\mapsto\varrho(x)\alpha$, and since
$w.\alpha=\varepsilon_n(w)\alpha$ for all $w\in S_n$, this
gives an isomorphism
\beq \label{linkeqn}
\varepsilon_n\otimes\CalV{(1,n)}\cong A^{n-1}(M(n))
\eeq
of representations of $S_n$.
Under this isomorphism, the $r=1$ case of Corollary \ref{newbasiscor}
corresponds to \cite[Proposition (3.3)]{lehrersolomon}, and the main result
\cite[Theorem (3.9)]{lehrersolomon} becomes the following fact about
$\CalV{(1,n)}$. Let $\psi_n:\mu_n\hookrightarrow\C^\times$ be the inclusion 
map, and simultaneously identify 
$\mu_n$ with a subgroup of $S_n$, by sending a generator $\zeta_0\in\mu_n$
to the $n$-cycle $(1,2,\cdots,n)$.
\bth \label{lsthm}
As a representation of $S_n$, $\CalV{(1,n)}\cong\Ind_{\mu_n}^{S_n}(\psi_n)$.
\eth
\bpf
To emphasize that Lehrer's and Solomon's proof is purely combinatorial,
we will sketch its translation into the context of $\CalV{(1,n)}$.
Since $\dim\CalV{(1,n)}=(n-1)!=|S_n/\mu_n|$, it suffices to show that
$\CalV{(1,n)}$ is generated as a $\C S_n$-module by some
$x$ satisfying $\theta.x=\psi_n(\theta)x$ for all $\theta\in\mu_n$.
An element which indeed has the latter property is 
$x=(\sum_{p=0}^{n-1}\zeta_0^{-p}(1,2,\cdots,n)^p).[T_0]$.
The claim is that
\beq \label{lseqn}
x=(1-c_{n-1}\zeta_0^{-1})(1-c_{n-2}\zeta_0^{-1})
\cdots(1-c_2\zeta_0^{-1})(1-c_1\zeta_0^{-1}).[T_0],
\eeq
where $c_i$ is the $i$-cycle $(1,2,\cdots,i)$.
Since each $1-c_i\zeta_0^{-1}$ for $i\leq n-1$ is an invertible
element of $\C S_n$, and $\CalV{(1,n)}$ is certainly generated 
as a $\C S_n$-module by $[T_0]$ (by Corollary \ref{newbasiscor}),
this finishes the proof. To prove the claim \eqref{lseqn}, 
we define elements $b_{n,p}\in\C S_{n-1}$ by the rule
\beq
(1-c_{n-1}t)(1-c_{n-2}t)
\cdots(1-c_2 t)(1-c_1 t)=\sum_{p=0}^{n-1}b_{n,p}t^p.
\eeq
It suffices to show that $c_n^p.[T_0]=b_{n,p}.[T_0]$, for all
$0\leq p\leq n-1$. The $p=0$ case is trivial, and the $p=n-1$
case is equivalent to $w_0.[T_0]=(-1)^{n-1}[T_0]$ where
$w_0$ is the permutation $i\mapsto n+1-i$, which is an immediate consequence
of \eqref{secondtreereleqn}. 
For the remaining cases $1\leq p\leq n-2$ we use induction on $n$, 
and the obvious recurrence
$b_{n,p}=b_{n-1,p}-c_{n-1}b_{n-1,p-1}$.
By the induction hypothesis, we have 
\[
b_{n-1,p}.[1\to 2\to\cdots\to n-1]=[p+1\to\cdots\to n-1\to 
1\to\cdots\to p],
\]
which, bearing in mind $b_{n-1,p}\in\C S_{n-2}$, implies that
\bes
\begin{split}
b_{n-1,p}.[T_0]&=-b_{n-1,p}.\left[\begin{array}{ccccc}
1&\to&2&\to\cdots\to&n-1\\
&&&&\uparrow\\
&&&&n\end{array}\right]\\
&=-\left[\begin{array}{ccccccc}
p+1&\to\cdots\to&n-1&\to&1&\to\cdots\to&p\\
&&\uparrow&&&&\\
&& n &&&&\end{array}\right].
\end{split}
\ees
Hence also
\bes
c_{n-1}b_{n-1,p-1}.[T_0]=
-\left[\begin{array}{ccccccc}
p+1&\to\cdots\to&n-1&\to&1&\to\cdots\to&p\\
&&&&\uparrow&&\\
&&&& n &&\end{array}\right].
\ees
The desired equation 
$c_n^p.[T_0]=b_{n-1,p}.[T_0]-c_{n-1}b_{n-1,p-1}.[T_0]$
is now a case of \eqref{firsttreereleqn}.
\epf
Combining our embeddings $\mu_n\hookrightarrow S_n$ and
$\mu_r\times S_n\hookrightarrow W(r,n)$, we can embed $\mu_r\times\mu_n$
in $W(r,n)$.
\bcr \label{combcor}
As a representation of $W(r,n)$,
\[ \CalV{(r,n)}\cong\Ind_{\mu_r\times\mu_n}^{W(r,n)}(1\times
\psi_n). \]
\ecr
\bpf
This follows from (5) of Proposition \ref{indprop} and Theorem \ref{lsthm}.
\epf

We will need a generalization of Theorem \ref{lsthm}. Suppose for the
remainder of this section that $r\mid n$, and embed $W(r,n/r)$ 
in $S_{n}$ as the centralizer of the product of $n/r$ disjoint $r$-cycles.
Hence $\mu_r\times\mu_{n/r}$ is also embedded in $S_n$.
\bth \label{mysterythm}
As a representation of $S_{n}$, 
\[ \CalV{(1,n)}\cong
\Ind_{\mu_r\times\mu_{n/r}}^{S_{n}}(\psi_r^{n/r}\times\psi_{n/r})
\cong\Ind_{W(r,n/r)}^{S_{n}}(\rmprod_{n/r}\otimes
\CalV{(r,n/r)}). \]
\eth
\bpf
The second isomorphism follows from Corollary \ref{combcor}.
It would be good to realize the first isomorphism explicitly, as
we did in the $r=1$ case (Theorem \ref{lsthm}), by finding a generator $x$
for $\CalV{(1,n)}$ which satisfies $\zeta.x=\psi_r(\zeta)^{n/r}x$
for all $\zeta\in\mu_r$ and $\theta.x=\psi_{n/r}(\theta)$ for all
$\theta\in\mu_{n/r}$. However, I have not been able to do this. 
We will instead prove that the characters coincide:
\beq \label{mysteryeqn}
\chi(\Ind_{\mu_{n}}^{S_{n}}(\psi_{n}))=
\chi(\Ind_{\mu_r\times\mu_{n/r}}^{S_{n}}(\psi_r^{n/r}\times\psi_{n/r})).
\eeq
Both sides are clearly supported on the elements of 
cycle type $(d^{n/d})$ where
$d\mid n$. The value of the left-hand side
at such an element, when multiplied by the index of the centralizer, 
is $\sum_{\zeta\in\mu_{d}^\circ}\zeta=
\mu(d)$, where $\mu_d^\circ$ denotes the set of primitive
$d$th roots of $1$, and $\mu(d)$ is the M\"obius $\mu$-function.
The value of the right-hand side at such an element,
when multiplied by the index of the centralizer, is
\bes
\begin{split}
\sum_{\substack{e\mid r\\f\mid(n/r)\\\lcm(e,f)=d}}
\sum_{\substack{\eta\in\mu_e^\circ\\\theta\in\mu_{f}^\circ}}\eta^{n/r}\theta
&=\sum_{\substack{e\mid r\\f\mid(n/r)\\\lcm(e,f)=d}}\mu(f)
\sum_{\eta\in\mu_e^\circ}\eta^{n/r}\\
&=\sum_{\substack{e\mid r\\f\mid(n/r)\\\lcm(e,f)=d}}
\frac{\mu(f)\mu(\frac{e}{\gcd(n/r,e)})\phi(e)}{\phi(\frac{e}{\gcd(n/r,e)})},
\end{split}
\ees
where $\phi$ is Euler's function. Since both formulas are multiplicative,
it suffices to prove that they are equal when $r=p^a$, $n/r=p^b$, $d=p^c$
for some prime $p$ and $a,b,c\in\N$, $c\leq a+b$. That is, we must prove
\beq 
\mu(p^c)=\sum_{\substack{0\leq x\leq a\\0\leq y\leq b\\\max\{x,y\}=c}}
\frac{\mu(p^y)\mu(p^{x-\min\{b,x\}})\phi(p^x)}{\phi(p^{x-\min\{b,x\}})}.
\eeq
If $c=0$, both sides are clearly $1$. If $c=1$ and $a=0$, the only term
on the right-hand side is $x=0$, $y=1$, and both sides are $-1$.
If $c=1$ and $b=0$, the only term on the right-hand side is $x=1$, $y=0$,
and both sides are again $-1$. If $c=1$ and $a,b\geq 1$, there are three
terms on the right-hand side, those where $x,y\in\{0,1\}$ are not both $0$;
their sum is once again $-1$. So we may assume that $c\geq 2$, in which case
the left-hand side is $0$. Because of the $\mu$-function factors, 
the only nonzero terms on the right-hand side are those where $y\leq 1$
and $x\leq b+1$. The first of these conditions forces $x=c$, so if $c>b+1$
or $c>a$ there are no nonzero terms at all, and we are done.
The only remaining case is that $c\geq 2$, $c\leq a$, $c\leq b+1$, 
in which case
the two nonzero terms $x=c$, $y=0$ and $x=c$, $y=1$ cancel each other out. 
\epf
\section{Cohomology of the hyperplane complement $T(r,n)$}
As in the previous section, $r$ and $n$ denote positive integers. Define
\[  T(r,n):= \{(z_1,z_2,\cdots,z_n)\in \C^{n}\,|\,
z_i\neq 0,\,\forall i,\, z_i^r\neq z_j^r,\,\forall i\neq j\}. \]
This is the complement of the hyperplanes
$z_i=0$ for $1\leq i\leq n$ and $z_i=\eta z_j$ for $1\leq i\neq j\leq n$,
$\eta\in\mu_r$. The wreath product $W(r,n)$ acts on $T(r,n)$ by the
restriction of its standard representation on $\C^n$. Indeed, if $r\geq 2$,
the hyperplanes we have removed are exactly the reflecting hyperplanes for 
the image $G(r,1,n)$ of this representation; if $r=1$, we have also removed 
the hyperplanes $z_i=0$.

Let $\Omega^\bullet(T(r,n))$ be the graded anti-commutative algebra of
regular differential forms on $T(r,n)$. Thus $\Omega^0(T(r,n))=\CalO(T(r,n))$
is the algebra of regular functions on $T(r,n)$, or in other words the subring
$\C[z_i^{\pm 1}, (z_i-\eta z_j)^{-1}]$ of $\C(z_1,\cdots,z_n)$, and for
$1\leq p\leq n$, $\Omega^p(T(r,n))$ is the free $\Omega^0(T(r,n))$-module with
basis $dz_{i_1}\wedge\cdots\wedge dz_{i_p}$, for $1\leq i_1<\cdots<i_p\leq n$.
We have the usual differential 
$d:\Omega^\bullet(T(r,n))\to\Omega^\bullet(T(r,n))$
raising degrees by $1$. As with any nonsingular affine complex variety,
the cohomology of the cochain complex
$(\Omega^\bullet(T(r,n)),d)$ is merely the cohomology $H^p(T(r,n),\C)$.
The natural action of $W(r,n)$ on $\Omega^\bullet(T(r,n))$ preserves 
the grading,
the wedge product, and the differential. Explicitly, the action
on $\Omega^0(T(r,n))$ is given by the formulas:
\beq
\begin{split} 
\unz.f(z_1,\cdots,z_n)
&=f(\zeta_1^{-1}z_1,\cdots,\zeta_n^{-1}z_n),\\
w.f(z_1,\cdots,z_n)
&=f(z_{w(1)},\cdots,z_{w(n)}),
\end{split}
\eeq
and the action on $\Omega^p(T(r,n))$ uses the additional rule
$w.dz_i=dz_{w(i)}$.

As mentioned in the introduction,
the cohomology of a hyperplane
complement can be described quite explicitly
by results of Brieskorn, Orlik and Solomon; in our case, this description is
as follows. Let $A^\bullet(T(r,n))$
be the subalgebra of
$\Omega^\bullet(T(r,n))$ generated by the following $1$-forms:
\[ \omega_i:= \frac{dz_i}{z_i},\,\forall i,\text{ and }
\omega_{i,j,\eta}:=\frac{dz_i-\eta dz_j}{z_i-\eta z_j},\,\forall i\neq j,\
\eta\in\mu_r. \]
Since these generators are closed, $d\alpha=0$ for all 
$\alpha\in A^\bullet(T(r,n))$.
\bth \label{brieskornthm} 
\textup{(Brieskorn, see \cite[Theorem 5.89]{orlikterao}.)}
The inclusion $(A^\bullet(T(r,n)),0)\hookrightarrow
(\Omega^\bullet(T(r,n)),d)$ of cochain complexes is a quasi-isomorphism.
In other words,
$A^p(T(r,n))\cong H^p(T(r,n),\C)$ via
the map sending a form $\alpha$ to its cohomology class.
\eth
\bth \label{nbcthm}
\textup{(Orlik-Solomon, see 
\cite[Theorems 3.43 and 3.126]{orlikterao}.)}
$A^\bullet(T(r,n))$ has a basis consisting of all products
$\alpha_1\wedge\alpha_2\wedge\cdots\wedge\alpha_n$, where each $\alpha_i$
is either $1$, $\omega_i$, or $\omega_{i,j,\eta}$ for some $j>i$ and
$\eta\in\mu_r$.
\eth
\noindent
To obtain the basis given here as the `no-broken-circuit' basis
of \cite[Theorem 3.43]{orlikterao}, order the hyperplanes so that
the hyperplanes $\{z_i=\eta z_j\}$ for $j>i$ come first, in lexicographic
order of $(i,j)$, and then the hyperplanes $\{z_i=0\}$, in order of $i$.
If $r=1$, we recover the basis used by Lehrer in \cite{lehrerone}
(identify $A^\bullet(T(1,n))$ with the Orlik-Solomon algebra of $S_{n+1}$
by rewriting $\omega_i$ as $\omega_{i,n+1}$).
If $r=2$, we recover the basis used by Lehrer in \cite{hyperoctahedral}.

Of course, the basis in Theorem \ref{nbcthm} is not $W(r,n)$-stable,
so it is necessary to consider the linear relations satisfied by
the set of all $W(r,n)$-translates of basis elements. This set can
be parametrized by an enhanced version of the trees considered in the
previous section.
\bdf
A \textbf{directed forest} on $\{1,\cdots,n\}$ is a directed graph with vertex
set $\{1,\cdots,n\}$ such that each connected component is a directed tree
on the vertices it contains. A \textbf{decorated directed forest} $F$
on $\{1,\cdots,n\}$ is a directed forest on $\{1,\cdots,n\}$ equipped with,
firstly, a labelling of each edge by an element of $\mu_r$, and, secondly, 
a partition of the set of roots into \textbf{open} and \textbf{closed} roots. 
We write $i\overset{\eta}{\rightarrow}j$ (with a subscript $F$ when necessary)
to mean that
there is an edge from $i$ to $j$ in $F$ whose label is $\eta\in\mu_r$.
We say that $F$ is \textbf{rectified} if $i\overset{\zeta}{\rightarrow} j$
implies $i<j$.
Let $\CalF(r,n)$ be the set of decorated directed forests on $\{1,\cdots,n\}$,
and let $\CalF(r,n)^{k,l}$ be the subset of $\CalF(r,n)$ consisting of
those forests with $k$ edges and $l$ closed roots (and hence $n-k-l$
open roots).
Write $\CalF(r,n)^\circ$ and $\CalF(r,n)^{\circ,k,l}$ for the 
rectified subsets of these.
\edf
We define an action of $W(r,n)$ on $\CalF(r,n)$ by the same rules as for
$\CalT(r,n)$. (As well as permuting the vertices, $S_n$ respects the
openness of roots; in other words, $i$ is an open root of $F$
if and only if $w(i)$ is an open root of $w.F$.)
Clearly each $\CalF(r,n)^{k,l}$ is a $W(r,n)$-stable subset.

We can now define a suitably signed version of the set of $W(r,n)$-translates
of the Orlik-Solomon basis elements.
\bdf
For any decorated directed forest $F\in\CalF(r,n)^{k,l}$,
and any $i\in\{1,\cdots,n\}$, define the sign
\[ \varepsilon_i(F):=
(-1)^{i-1-|\{\text{open roots of $F$ which are $<i$}\}|}. \]
Define $\alpha(F)\in A^{k+l}(T(r,n))$ by
$\alpha_1(F)\wedge\alpha_2(F)\wedge\cdots\wedge\alpha_n(F)$, where
\[ \alpha_i(F):=\left\{\begin{array}{cl}
\varepsilon_i(F),&\text{ if $i$ is an open root of $F$,}\\
\omega_i,&\text{ if $i$ is a closed root of $F$,}\\
\omega_{i,j,\eta},&\text{ if $i\overset{\eta}{\rightarrow}j$ in $F$.}
\end{array}\right. \]
\edf
Note that if $i$ is not an open root of $F$, $\varepsilon_i(F)$ is exactly
the sign incurred in moving the $\alpha_i(F)$ factor of the wedge product to 
the front.

Since the wedge products in Theorem \ref{nbcthm} are, up to sign, those
$\alpha(F)$ for $F$ rectified, we can rewrite that result as follows.
\bth \label{mynbcthm}
$A^p(T(r,n))$ has basis $\{\alpha(F)\,|\,F\in\bigcup_{k+l=p}
\CalF(r,n)^{\circ,k,l}\}$.
\eth
\bex
The twelve elements of $\CalF(2,2)$, with the corresponding elements
of $A^\bullet(2,2)$, are as follows; asterisks
indicate closed roots. (Note that $2\overset{\pm 1}{\longrightarrow}1$
and $2\overset{\pm 1}{\longrightarrow}1^*$ are not rectified.)
\[ \begin{array}{|c|c|}
\hline
1\qquad 2&1\\
\hline
1\overset{\pm 1}{\longrightarrow}2&-\omega_{1,2,\pm 1}\\
\hline
2\overset{\pm 1}{\longrightarrow}1&\omega_{2,1,\pm 1}\\
\hline
1\qquad 2^* & \omega_2\\
\hline
1^*\qquad 2 & -\omega_1\\
\hline
1\overset{\pm 1}{\longrightarrow}2^*&\omega_{1,2,\pm 1}\wedge\omega_2\\
\hline
2\overset{\pm 1}{\longrightarrow}1^*&\omega_1\wedge\omega_{2,1,\pm 1}\\
\hline
1^*\qquad 2^* & \omega_1\wedge\omega_2\\
\hline
\end{array}
\]
\eex

The justification for the above definition of $\alpha(F)$
is the following result, in which the sign 
$\varepsilon(w,F)$ is defined by
\[ \varepsilon(w,F):=(-1)^{|\{(i,j)\,|\,i,j\text{ open roots of }F,\ i<j,\ 
w(i)>w(j)\}|}. \]
\bpr \label{actionprop}
For $\unz\in\mu_r^n$, $w\in S_n$ and $F\in\CalF(r,n)$,
\[ \alpha(\unz w.F)=\varepsilon_n(w)\,\varepsilon(w,F)\,
\unz w.\alpha(F). \]
\epr
\bpf
That $\unz.\alpha(F)=\alpha(\unz.F)$
is obvious from the fact that
$\unz.\omega_{i,j,\eta}=
\omega_{i,j,\eta\zeta_i\zeta_j^{-1}}$.
So we need only show that $\alpha(w.F)=\varepsilon_n(w)\varepsilon(w,F)
w.\alpha(F)$,
which it suffices to prove in the case that $w$ is a simple transposition
$s_i=(i,i+1)$, since $\varepsilon(ww',F)=\varepsilon(w,w'F)\varepsilon(w',F)$.
It is clear that $\alpha_j(s_i.F)=s_i.\alpha_j(F)$
for all $j\neq i,i+1$, so we are reduced to showing
\[ \alpha_i(s_i.F)\wedge\alpha_{i+1}(s_i.F)=\left\{\begin{array}{cl}
s_i.\alpha_i(F)\wedge s_i.\alpha_{i+1}(F),&\text{ if $i$, $i+1$
are open}\\
&\qquad\text{ roots of $F$}\\
-s_i.\alpha_i(F)\wedge s_i.\alpha_{i+1}(F),&\text{ otherwise.}
\end{array}\right. \]
Now if $i$ and $i+1$ are open roots of $F$,
both sides are clearly $1$. If
neither $i$ nor $i+1$ is an open root of $F$,
then $\alpha_i(s_i.F)=s_i.\alpha_{i+1}(F)$ and
$\alpha_{i+1}(s_i.F)=s_i.\alpha_i(F)$, both being $1$-forms, so the claim
holds. If $i+1$ is an open root of $F$ and $i$ is not, then
$\alpha_{i+1}(s_i.F)=s_i.\alpha_{i}(F)$ but
\bes
\begin{split}
\varepsilon_{i}(s_i.F)
&=(-1)^{i-1-|\{\text{open roots of $s_i.F$ which are $<i$}\}|}\\
&=(-1)^{i-1-|\{\text{open roots of $F$ which are $<i+1$}\}|}
=-\varepsilon_{i+1}(F),
\end{split}
\ees
so again the claim holds. The case where $i$ is an open root
of $F$ and $i+1$ is not is similar. 
\epf

We now describe the linear relations
satisfied by the elements $\alpha(F)$,
when $F$ is not necessarily rectified; these were the motivation
for the relations used in the previous section. 
\bpr \label{relprop}
The following hold in $A^\bullet(T(r,n))$.
\ben
\item $\alpha(F_1)+\alpha(F_2)=\alpha(F_3)$,
if $F_1,F_2,F_3\in\CalF(r,n)$
are identical except for the edges between three vertices $i$, $j$, and $k$,
where we have

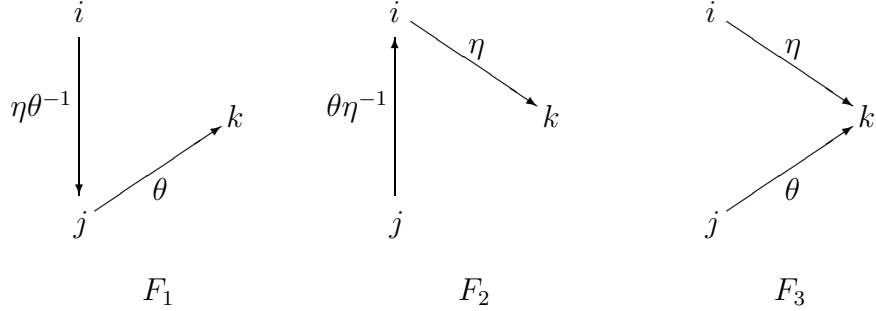
\setlength{\unitlength}{0.7mm}
\begin{picture}(160,70)(10,0)
\put(9,58){$i$}
\put(9,18){$j$}
\put(38,38){$k$}
\put(10,55){\vector(0,-1){30}}
\put(13,22){\vector(3,2){24}}
\put(-3,40){$\eta\theta^{-1}$}
\put(24,24){$\theta$}
\put(22,5){$F_1$}
\put(69,58){$i$}
\put(69,18){$j$}
\put(98,38){$k$}
\put(70,25){\vector(0,1){30}}
\put(73,58){\vector(3,-2){24}}
\put(57,40){$\theta\eta^{-1}$}
\put(84,52){$\eta$}
\put(82,5){$F_2$}
\put(129,58){$i$}
\put(129,18){$j$}
\put(158,38){$k$}
\put(133,58){\vector(3,-2){24}}
\put(133,22){\vector(3,2){24}}
\put(144,24){$\theta$}
\put(144,52){$\eta$}
\put(142,5){$F_3$}
\end{picture}

\noindent
for some $\eta,\theta\in\mu_r$. (The assertion that they are otherwise
identical includes the fact that if $k$ is a root of any, hence all, of these
forests, it is either open in all or closed in all.)
\item $\alpha(F)+\alpha(F')=0$,
if $F,F'\in\CalT(r,n)$
are identical except that $j$ is an open root of $F$,
$i$ is an open root of $F'$, and we have $i\overset{\eta}{\rightarrow} j$
in $F$, $j\overset{\eta^{-1}}{\rightarrow} i$ in $F'$ for some $\eta\in\mu_r$.
\item $\alpha(F)+\alpha(F')=\alpha(F'')$,
if $F,F',F''\in\CalF(r,n)$ are
identical except that $j$ is a closed root of $F$,
$i$ is a closed root of $F'$, both are closed roots of $F''$,
and we have $i\overset{\eta}{\rightarrow} j$
in $F$, $j\overset{\eta^{-1}}{\rightarrow} i$ in $F'$ for some $\eta\in\mu_r$
(there is of course no edge between $i$ and $j$
in $F''$).
\een
\epr
\bpf
In (1) and (3),
the factors in various positions of the wedge product
are forms of the same degrees in all three terms, so we may reorder
the factors in the same way in all three terms so as to bring the factors
in which the terms differ to the front, and then neglect the other factors.
Therefore to prove these relations we need only check:
\beq
\begin{split}
\omega_{i,j,\eta\theta^{-1}}\wedge\omega_{j,k,\theta}
+\omega_{i,k,\eta}\wedge\omega_{j,i,\theta\eta^{-1}}
&=\omega_{i,k,\eta}\wedge\omega_{j,k,\theta},\text{ and}\\
\omega_{i,j,\eta}\wedge\omega_j+\omega_i\wedge\omega_{j,i,\eta^{-1}}
&=\omega_i\wedge\omega_j.
\end{split}
\eeq
These are known relations in $A^\bullet(T(r,n))$ (besides being trivial
to verify directly).
In (2), reorder
the factors of the wedge products so that the factor corresponding to
$i$ comes first, that corresponding to $j$ comes second, and the others
(which are the same in both terms) follow in order. In the case of
$\alpha(F)$, the sign incurred by this reordering is 
$\varepsilon_i(F)$, and the reordered
product is $\omega_{i,j,\eta}\wedge\varepsilon_j(F)\wedge\cdots$.
In the case of $\alpha(F')$, the sign incurred is $\varepsilon_j(F')$,
and the reordered product is
$\varepsilon_i(F')\wedge\omega_{j,i,\eta^{-1}}\wedge\cdots$. Since
$\omega_{i,j,\eta}=\omega_{j,i,\eta^{-1}}$, we need only show that
$\varepsilon_i(F)\varepsilon_j(F)+\varepsilon_i(F')\varepsilon_j(F')=0$.
But if $i<j$, then $\varepsilon_i(F')=\varepsilon_i(F)$ while
$\varepsilon_j(F')=-\varepsilon_j(F)$, because $i$ is an open root
of $F'$ but not of $F$; the $j<i$ case is similar.
\epf

From the viewpoint of the previous section,
relation (3) of Proposition \ref{relprop}
appears to be of the wrong form. However, we can
make $\alpha(F)+\alpha(F')$ zero as expected if we work 
modulo forests with a greater number of closed roots (such as $F''$).
So we define for each $p$ a $W(r,n)$-stable filtration 
\[ A^p(T(r,n))=A^p(T(r,n))_0\supseteq A^p(T(r,n))_1\supseteq A^p(T(r,n))_2
\cdots, \]
where
\[ A^p(T(r,n))_m:=\mathrm{span}\{\alpha(F)\,|\,F\in
\bigcup_{\substack{k+l=p\\l\geq m}}\CalF(r,n)^{k,l}\}. \]
(In the $r=2$ case, a similar filtration was used in \cite{douglass}.)
Also define $A^{k,l}(T(r,n)):=
A^{k+l}(T(r,n))_l/A^{k+l}(T(r,n))_{l+1}$, and for $F\in\CalF(r,n)^{k,l}$,
let $[F]$ denote the image of $\alpha(F)$
in $A^{k,l}(T(r,n))$. We can now state the main result of this section.
\bth \label{relthm} Let $r,n,p,k,l$ be as above.
\ben
\item $A^p(T(r,n))$ can be defined abstractly as the vector space spanned
by $\{\alpha(F)\,|\,F\in\bigcup_{k+l=p}\CalF(r,n)^{k,l}\}$ subject to the 
relations in Proposition \ref{relprop}.
\item $A^{k,l}(T(r,n))$ can be defined abstractly 
as the vector space spanned by
$\{[F]\,|\,F\in\CalF(r,n)^{k,l}\}$ subject to the following relations:
\begin{itemize}
\item $[F_1]+[F_2]=[F_3]$ whenever $F_1,F_2,F_3\in\CalF(r,n)^{k,l}$ are as in
(1) of Proposition \ref{relprop};
\item $[F]+[F']=0$ whenever 
$F,F'\in\CalF(r,n)^{k,l}$ are as in (2) or (3) of Proposition \ref{relprop}.
\end{itemize}
\item $A^{k,l}(T(r,n))$ has basis $\{[F]\,|\,F\in\CalF(r,n)^{\circ,k,l}\}$.
\item For $\unz\in\mu_r^n$, $w\in S_n$ and $F\in\CalF(r,n)^{k,l}$,
\[ \unz w.[F]=\varepsilon_n(w)\,\varepsilon(w,F)\,[\unz w.F]. \]
\item  As a representation of $W(r,n)$,
\[ A^p(T(r,n))\cong\bigoplus_{k+l=p}A^{k,l}(T(r,n)). \]
\een
\eth
\bpf
The vector space defined abstractly as in (1) certainly maps to
$A^p(T(r,n))$, since the stated relations do hold. So 
to prove (1), it suffices to show that the basis elements
$\{\alpha(F)\,|\,F\in\bigcup_{k+l=p}\CalF(r,n)^{\circ,k,l}\}$
of Theorem \ref{mynbcthm}
span the abstract vector space also. 
This follows from the same sort of rectification
procedure as in Lemma \ref{spanlemma}: for $F\in\CalF(r,n)^{k,l}$,
we use relations (2) and (3) of Proposition \ref{relprop} to shift the root
of each tree to its largest vertex, and any forests which arise from
the right-hand side of relation (3) can be neglected,
by downward induction on the number of closed roots. Once the root
of each tree is its largest vertex, the rectification proceeds using
relation (1) as in Lemma \ref{spanlemma}. So part (1) of the present result
is proved. Moreover, this procedure expresses each $\alpha(F)$
as a linear combination of $\alpha(F')$'s for $F'$ rectified with at least
as many closed roots as $F$. Therefore
$A^p(T(r,n))_m$ has basis $\{\alpha(F)\,|\,
F\in\bigcup_{\substack{k+l=p\\l\geq m}}\CalF(r,n)^{\circ,k,l}\}$, from which
part (3) follows immediately. 
Part (2) follows from (3) by
the same argument as for (1). Part (4) clearly follows from 
Proposition \ref{actionprop}, and part (5) from complete reducibility
of finite-dimensional representations of $W(r,n)$.
\epf

We can now relate $A^{k,l}(T(r,n))$ to the representations considered in the
previous section. Recall that
if $\lambda=(\lambda_1,\cdots,\lambda_\ell)$ is a 
partition in the usual combinatorial sense,
the stabilizer in $S_{|\lambda|}$ of a set partition of 
$\{1,\cdots,|\lambda|\}$ into parts of sizes $\lambda_1$, $\cdots$, 
$\lambda_\ell$
is isomorphic to 
\[ (S_{\lambda_1}\times\cdots\times S_{\lambda_\ell})\rtimes
(S_{m_1(\lambda)}\times S_{m_2(\lambda)}\times\cdots), \] 
where $S_{m_i(\lambda)}$ acts by permuting the $S_{\lambda_a}$ factors where 
$\lambda_a=i$.
The subgroups of $W(r,n)$ referred to in the following result are defined
similarly.
\bcr \label{bigindcor}
As a representation of $W(r,n)$, $\varepsilon_n\otimes A^{k,l}(T(r,n))$
is isomorphic to the following direct sum:
\[ \bigoplus_{\substack{\lambda^1,\lambda^2\\|\lambda^1|+|\lambda^2|=n
\\\ell(\lambda^1)=n-k-l\\
\ell(\lambda^2)=l}}
\negthickspace\negthickspace
\Ind_{\substack{((W(r,\lambda_1^1)\times\cdots\times W(r,\lambda_{n-k-l}^1))\\
\,\rtimes(S_{m_1(\lambda^1)}\times S_{m_2(\lambda^1)}\times\cdots))\\
\negthickspace
\times((W(r,\lambda_1^2)\times\cdots\times W(r,\lambda_{l}^2))\\
\,\rtimes(S_{m_1(\lambda^2)}\times S_{m_2(\lambda^2)}\times\cdots))}}^{W(r,n)}
\left(\negthickspace\begin{array}{c}
\varepsilon\otimes
\CalV{(r,\lambda_1^1)}\otimes\cdots\otimes \CalV{(r,\lambda_{n-k-l}^1)}\\
\phantom{\varepsilon}\otimes
\CalV{(r,\lambda_1^2)}\otimes\cdots\otimes \CalV{(r,\lambda_{l}^2)}
\end{array}\negthickspace\right), \]
where $W(r,\lambda_a^j)$ acts on the $\CalV{(r,\lambda_a^j)}$
factor, $S_{m_i(\lambda^j)}$ acts by permuting the $\CalV{(r,\lambda_a^j)}$
factors where $\lambda_a^j=i$, and $\varepsilon$ denotes
the product of the sign characters of the $S_{m_i(\lambda^1)}$ components.
\ecr
\bpf
For any partition of the set $\{1,\cdots,n\}$ into nonempty subsets
$B_1,\cdots,B_{n-k-l}$, $C_1,\cdots,C_l$, let 
$\CalF(r,n)^{k,l}_{(\{B_i\},\{C_j\})}\subset\CalF(r,n)^{k,l}$
be the subset consisting of forests in which the trees with open roots
have vertex sets $B_1,\cdots,B_{n-k-l}$
and the trees with closed roots have vertex sets $C_1,\cdots,C_l$.
Let $A^{k,l}(T(r,n))_{(\{B_i\},\{C_j\})}$
be the subspace of $A^{k,l}(T(r,n))$ spanned by the elements
$[F]$ for $F\in\CalF(r,n)^{k,l}_{(\{B_i\},\{C_j\})}$.
Since every relation in (2) of Theorem \ref{relthm} involves forests in
the same subset $\CalF(r,n)^{k,l}_{(\{B_i\},\{C_j\})}$, 
we have a direct sum decomposition
\beq 
A^{k,l}(T(r,n))=\bigoplus_{(\{B_i\},\{C_j\})} 
A^{k,l}(T(r,n))_{(\{B_i\},\{C_j\})}.
\eeq
Comparing the relations in Theorem \ref{relthm}
to those in the definition of $\CalV{(r,n)}$, we see that there 
is an isomorphism of vector spaces
\beq \label{splittingeqn}
\begin{split}
A^{k,l}(T(r,n))_{(\{B_i\},\{C_j\})}&\isomto
\CalV{(r,B_1)}\otimes\cdots\otimes \CalV{(r,B_{n-k-l})}\\
&\qquad\qquad\otimes
\CalV{(r,C_1)}\otimes\cdots\otimes \CalV{(r,C_l)}
\end{split}
\eeq
sending $[F]$ for $F\in\CalF(r,n)^{k,l}_{(\{B_i\},\{C_j\})}$ to 
\[ [F|_{B_1}]\otimes\cdots\otimes[F|_{B_{n-k-l}}]\otimes
[F|_{C_1}]\otimes\cdots\otimes[F|_{C_l}]. \]
Here $\CalV{(r,X)}$ is defined in the same way as
$\CalV{(r,n)}$ but with vertex set $X$, and $F|_{X}$ means
the tree in the forest $F$ whose vertex set is $X$.
Now the stabilizer of $(\{B_i\},\{C_j\})$ in $W(r,n)$ is a subgroup
of the type mentioned in the statement, where $\lambda^1$ is formed
from the sizes of the $B_i$'s and $\lambda^2$ from the sizes of the $C_j$'s.
Thanks to (4) of Theorem \ref{relthm}, the action of this stabilizer
on $\varepsilon_n\otimes A^{k,l}(T(r,n))_{(\{B_i\},\{C_j\})}$
corresponds to the obvious action on the right-hand side of 
\eqref{splittingeqn}, together with a sign on the permutation of the open-root
factors. The result follows.
\epf
\bex
In the case $r=n=2$, Corollary \ref{bigindcor} gives the following isomorphisms
of representations of $W(2,2)$:
\bes
\begin{split}
A^{0,0}(T(2,2))&\cong\varepsilon_2\otimes 
\Ind_{(W(2,1)\times W(2,1))\rtimes S_2}^{W(2,2)}(\varepsilon\otimes
\CalV{(2,1)}\otimes\CalV{(2,1)})\cong 1,\\
A^{1,0}(T(2,2))&\cong\varepsilon_2\otimes 
\Ind_{W(2,2)}^{W(2,2)}(\CalV{(2,2)})\cong 1\oplus\rmprod_2,\\
A^{0,1}(T(2,2))&\cong\varepsilon_2\otimes 
\Ind_{W(2,1)\times W(2,1)}^{W(2,2)}(\CalV{(2,1)}\otimes\CalV{(2,1)})
\cong 1\oplus\varepsilon_2,\\
A^{1,1}(T(2,2))&\cong\varepsilon_2\otimes 
\Ind_{W(2,2)}^{W(2,2)}(\CalV{(2,2)})\cong 1\oplus\rmprod_2,\\
A^{0,2}(T(2,2))&\cong\varepsilon_2\otimes 
\Ind_{(W(2,1)\times W(2,1))\rtimes S_2}^{W(2,2)}
(\CalV{(2,1)}\otimes\CalV{(2,1)})\cong\varepsilon_2.
\end{split}
\ees
Here we have used the obvious fact 
$\CalV(2,2)\cong\varepsilon_2\oplus\rmdet_2$.
\eex

In general, we can substitute Corollary \ref{combcor} 
into Corollary \ref{bigindcor} to obtain:
\bcr \label{mainindcor}
As a representation of $W(r,n)$, $\varepsilon_n\otimes A^{k,l}(T(r,n))$
is isomorphic to the following direct sum:
\[ \bigoplus_{\substack{\lambda^1,\lambda^2\\|\lambda^1|+|\lambda^2|=n\\
\ell(\lambda^1)=n-k-l\\\ell(\lambda^2)=l}}
\Ind_{\substack{(((\mu_r\times\mu_{\lambda_1^1})\times\cdots\times 
(\mu_r\times\mu_{\lambda_{n-k-l}^1}))
\rtimes(S_{m_1(\lambda^1)}\times S_{m_2(\lambda^1)}\times\cdots))\\
\!\times(((\mu_r\times\mu_{\lambda_1^2})\times\cdots
\times(\mu_r\times\mu_{\lambda_{l}^2}))
\rtimes(S_{m_1(\lambda^2)}\times S_{m_2(\lambda^2)}\times\cdots))}}^{W(r,n)}
(\varepsilon\psi),
\]
where $\psi$ is the character which takes the product of the 
$\mu_{\lambda_a^j}$ components, and $\varepsilon$ is the product of the sign
characters of the $S_{m_i(\lambda^1)}$ components.
\ecr
Summing over $k$ and $l$, we get an expression for
$A^\bullet(T(r,n))$ as a direct sum of
inductions of one-dimensional characters, which in the case $r=2$
is the same as \cite[(1.1)]{douglass}.

Another consequence of Corollary \ref{bigindcor} is a proof of the formula
for the character of $H^\bullet(T(r,n),\C)$ used in \cite{mytorus}.
Since there are already at least three proofs of this character 
formula in the literature (Hanlon's original proof of the equivalent
poset result in \cite{hanlonwreath}, Lehrer's proof in \cite{lehrertwo},
and the `equivariant inclusion-exclusion' argument of 
\cite[Theorem 9.4]{mywreath}), I will leave this to the reader's
imagination.
\section{Cohomology of $T(1,n)$ with coefficients in a local system}
Now let $a_1,\cdots,a_n\in\Q$,
and define the local system $\CalL_{\sum a_i\omega_i}$ on $T(1,n)$
as in the introduction. Let $r$ be a positive integer such that
$ra_i\in\Z$ for all $i$.
For any representation $V$ of $W(r,n)$, define
\bes
V_{(a_i)}:=\{v\in V\,|\,\unz.v=\zeta_1^{-ra_1}\cdots\zeta_n^{-ra_n}\, v,\,
\forall\unz\in\mu_r^n\}.
\ees
It is clear that $V_{(a_i)}$ is stable under $Z_{(a_i)}$, the subgroup
of $S_n$ which fixes $(a_i)$. If all $a_i$ equal $\frac{s}{r}$ for some
$s\in\Z$, we write $V_{(a_i)}$ as $V_\sr$; in this case 
$Z_{(a_i)}$ is $S_n$ itself.

We are interested in
$H^\bullet(T(r,n),\C)_{(a_i)}\cong A^\bullet(T(r,n))_{(a_i)}$, 
which as mentioned in the Introduction
is isomorphic to $H^\bullet(T(1,n),\CalL_{\sum a_i\omega_i})$ via the map
$\varphi:T(r,n)\to T(1,n):
(z_1,\cdots,z_n)\mapsto(z_1^r,\cdots,z_n^r)$.
To see this isomorphism explicitly, note that the pull-back
$\varphi^*:\Omega^\bullet(T(1,n))\to\Omega^\bullet(T(r,n))$ sends
$z_i$ to $z_i^r$, $dz_i$ to $rz_i^{r-1}dz_i$ and hence $\omega_i$ to
$r\omega_i$ for all $i$.
Recall that we have an isomorphism
\bes
H^p(T(1,n),\CalL_{\sum a_i\omega_i})\cong
\frac{\ker(d+\ell_{\sum a_i\omega_i}:\Omega^p(T(1,n))\to\Omega^{p+1}(T(1,n)))}
{\im(d+\ell_{\sum a_i\omega_i}:\Omega^{p-1}(T(1,n))\to\Omega^{p}(T(1,n)))},
\ees 
which is clearly $Z_{(a_i)}$-equivariant.
\bpr \label{cochainprop}
The map $\alpha\mapsto z_1^{ra_1}\cdots z_n^{ra_n}\varphi^*(\alpha)$
defines an isomorphism of complexes
$(\Omega^\bullet(T(1,n)),d+\ell_{\sum a_i\omega_i})\isomto
(\Omega^\bullet(T(r,n))_{(a_i)},d)$,
inducing $Z_{(a_i)}$-equivariant isoms 
$H^p(T(1,n),\CalL_{\sum a_i\omega_i})\isomto
H^p(T(r,n),\C)_{(a_i)}$.
\epr
\bpf
Since the action of $\mu_r^n$ on $\Omega^\bullet(T(r,n))$ respects $d$,
the second statement follows automatically from the first.
Now it is clear that $\varphi^*:\Omega^0(T(1,n))\to\Omega^0(T(r,n))$ is an 
isomorphism onto the $\mu_r^n$-invariant subspace
$\Omega^0(T(r,n))^{\mu_r^n}$. It follows that
$\varphi^*:\Omega^p(T(1,n))\isomto\Omega^p(T(r,n))^{\mu_r^n}$ for all
$p$. Hence we have an isomorphism of complexes
\beq \label{cochainisomeqn}
\varphi^*:(\Omega^\bullet(T(1,n)),d+\ell_{\sum a_i\omega_i})\isomto
(\Omega^\bullet(T(r,n))^{\mu_r^n},d+\ell_{\sum ra_i\omega_i}).
\eeq
Now for any $\alpha\in\Omega^p(T(r,n))$,
\bes
\begin{split}
z_1^{ra_1}\cdots z_n^{ra_n}(d&+\ell_{\sum ra_i\omega_i})(\alpha)\\
&=z_1^{ra_1}\cdots z_n^{ra_n}d\alpha+\sum_i z_1^{ra_1}\cdots
\widehat{z_i^{ra_i}}\cdots z_n^{ra_n}d(z_i^{ra_i})\wedge\alpha\\
&=d(z_1^{ra_1}\cdots z_n^{ra_n}\alpha).
\end{split}
\ees
Hence the map
$\alpha\mapsto z_1^{ra_1}\cdots z_n^{ra_n}\alpha$ gives an isomorphism
of complexes $(\Omega^\bullet(T(r,n)),d+\ell_{\sum ra_i\omega_i})
\isomto(\Omega^\bullet(T(r,n)),d)$, which clearly restricts to an
isomorphism $(\Omega^\bullet(T(r,n))^{\mu_r^n},d+\ell_{\sum ra_i\omega_i})
\isomto(\Omega^\bullet(T(r,n))_{(a_i)},d)$.
Combining this with the isomorphism \eqref{cochainisomeqn}, we have the result.
\epf

In the previous section we saw that the structure of $A^p(T(r,n))\cong
H^p(T(r,n),\C)$ was controlled by the forests in the set $\CalF(r,n)$.
The main point of this section is that the structure of $A^p(T(r,n))_{(a_i)}$
is controlled by the forests in the following subset of $\CalF(1,n)$.
\bdf
Let $\CalF(1,n;a_i)$ be the set of $F\in\CalF(1,n)$ such that
\beq \label{theintegraleqn}
\sum_{i\in T} a_i\in\Z,\text{ for all trees $T$ in $F$.}
\eeq
Define $\CalF(1,n;a_i)^{k,l}=\CalF(1,n;a_i)\cap\CalF(1,n)^{k,l}$,
and similarly $\CalF(1,n;a_i)^\circ$ and $\CalF(1,n;a_i)^{\circ,k,l}$.
\edf
Note that $\CalF(1,n;a_i)$ is nonempty if and only if
$\sum_i a_i\in\Z$. If all $a_i$ equal $\frac{s}{r}$, we write
$\CalF(1,n;a_i)$ as $\CalF(1,n;\frac{s}{r})$; it is then nonempty
if and only if $r\mid sn$, i.e.\ $\frac{r}{\gcd(s,r)}\mid n$. Moreover,
this being the case,
$\CalF(1,n;\frac{s}{r})^{k,l}$ is nonempty if and only if
$l\leq n-k\leq \frac{n\gcd(s,r)}{r}$, since the number of vertices in each tree
must be a multiple of $\frac{r}{\gcd(s,r)}$.

Now for any $F\in\CalF(1,n)^{k,l}$, we define an 
element of $A^{k+l}(T(r,n))_{(a_i)}$
in the most obvious way:
\beq
\beta(F):=\sum_{\unz\in\mu_r^n}
\zeta_1^{ra_1}\cdots\zeta_n^{ra_n}\unz.\alpha(F)=
\sum_{\unz\in\mu_r^n}\zeta_1^{ra_1}\cdots\zeta_n^{ra_n}\alpha(\unz.F).
\eeq
\bpr \label{betaprop}
\ben
\item $\beta(F)\neq 0\Longleftrightarrow F\in\CalF(1,n;a_i)$.
\item $A^\bullet(T(r,n))_{(a_i)}$ has basis 
$\{\beta(F)\,|\,F\in\CalF(1,n;a_i)^\circ\}$.
\item For $w\in Z_{(a_i)}$ and $F\in\CalF(1,n;a_i)$,
\[ w.\beta(F)=\varepsilon_n(w)\varepsilon(w,F)\beta(w.F). \]
\een
\epr
\bpf
We know that the set 
$\{\alpha(\unz.F)\,|\,\unz\in\mu_r^n\}$ is
linearly independent, although
each element can correspond to many different $\unz$.
(If $F$ is rectified, this is part of the basis in Theorem \ref{mynbcthm}, 
and the
linear independence in general follows by applying a suitable element 
of $S_n$.)
So it is clear that $\beta(F)\neq 0$ if and only if
\beq \label{nonvaneqn}
\sum_{\substack{\unz\in\mu_r^n\\\unz.F=F}}
\zeta_1^{ra_1}\cdots\zeta_n^{ra_n}\neq 0.
\eeq
But the subgroup of $\mu_r^n$ fixing $F$ consists of the elements
satisfying $\zeta_i=\zeta_j$ whenever $i$ and $j$ are vertices
of the same tree in $F$. So the left-hand side of \eqref{nonvaneqn}
factorizes into the product over the trees $T$ of
$\sum_{\zeta\in\mu_r}\zeta^{\sum_{i\in T}ra_i}$, which is nonzero
if and only if $r\mid \sum_{i\in T}ra_i$, i.e.\ $\sum_{i\in T}a_i\in\Z$. 
So (1) is proved.
Now from Theorem \ref{mynbcthm} and the fact that
any $F'\in\CalF(r,n)^{\circ}$ is of the form
$\unz.F$ for some $F\in\CalF(1,n)^{\circ}$, it is evident that
$A^{\bullet}(T(r,n))_{(a_i)}$ is spanned by 
$\{\beta(F)\,|\,F\in\CalF(1,n)^{\circ}\}$. By (1), this means that
$A^{\bullet}(T(r,n))_{(a_i)}$ is spanned by the nonzero elements
$\{\beta(F)\,|\,F\in\CalF(1,n;a_i)^{\circ}\}$. The latter elements
are clearly linearly independent, as they involve disjoint collections
of $\alpha(F')$'s, so (2) is proved. Part (3) is a consequence of
Proposition \ref{actionprop}.
\epf
As a corollary, we have 
\beq
\dim A^p(T(r,n))_{(a_i)}=\sum_{k+l=p}|\CalF(1,n;a_i)^{\circ,k,l}|.
\eeq
In particular, $A^p(T(r,n))_{(a_i)}=0$ unless $\sum_i a_i\in\Z$, and
\beq
A^p(T(r,n))_{(\frac{s}{n})}\neq 0 \Longleftrightarrow
n-\gcd(s,n)\leq p\leq n.
\eeq

As described in the Introduction,
we can write down explicit differential forms on $T(1,n)$
(not in the Orlik-Solomon algebra) which
give the basis of $H^p(T(1,n),\CalL_{\sum a_i\omega_i})$
corresponding under Proposition \ref{cochainprop}
to the basis
$\{\beta(F)\,|\,F\in\bigcup_{k+l=p}\CalF(1,n;a_i)^{\circ,k,l}\}$ of
$A^p(T(r,n))_{(a_i)}$.
For this we need some more notation associated with a forest
$F\in\CalF(1,n;a_i)$. Let $\simeq$ be the equivalence relation on
$\{1,\cdots,n\}$ generated by $\to$, so $i\simeq j$ means that $i$
and $j$ are vertices of the same tree.  
Let $\preceq$ be the partial order generated by $\to$,
so $i\preceq j$ means that there is a path
from $i$ to $j$
traversing each edge in the chosen direction.
Define
\bes
\begin{split}
b_j(F)&=-\lfloor\sum_{k\preceq j}a_k\rfloor
+\sum_{i\to j}\lceil\sum_{k\preceq i}a_k\rceil\\
&=-\sum_{k\simeq j}a_k+\lceil\sum_{\substack{k\simeq j\\k\not\preceq j}}a_k
\rceil+\sum_{i\to j}\lceil\sum_{k\preceq i}a_k\rceil.
\end{split}
\ees
(From the second expression it is clear that $b_j(F)$ is independent
of the orientation of edges in $F$.)
Then define
\bes 
\begin{split}
\bbeta(F)&:=
(\prod_{j=1}^n z_j^{b_j(F)})\,\bbeta_1(F)\wedge\cdots\wedge\bbeta_n(F)
\in\Omega^\bullet(T(1,n)),\text{ where}\\
\bbeta_i(F)&:=\left\{\begin{array}{cl}
\varepsilon_i(F) r,&\text{ if $i$ is an open root of $F$,}\\
\omega_i,&\text{ if $i$ is a closed root of $F$,}\\
\omega_{i,j},&\text{ if $i\to j$ in $F$ and $\sum_{k\preceq i}a_k\in\Z$,}\\
(z_i-z_j)^{-1}(\omega_i-\omega_j),&\text{ if $i\to j$ in $F$ 
and $\sum_{k\preceq i}a_k\not\in\Z$.}
\end{array}\right.
\end{split}
\ees
Note that $(z_i-z_j)^{-1}(\omega_i-\omega_j)=\frac{z_j dz_i-z_i dz_j}
{z_i z_j(z_i-z_j)}$, and that $\sum_{k\preceq i}a_k\in\Z$ if and only if
the edge $i\to j$ is \textbf{breakable}, in the sense that the new trees 
formed by deleting it still satisfy \eqref{theintegraleqn}.
\bpr \label{calcprop}
$\beta(F)=z_1^{ra_1}\cdots z_n^{ra_n}\,
\varphi^*(\bbeta(F)),\ \forall F\in\CalF(1,n;a_i)$.
\epr
\bpf
It is easy to see that
$\beta(F)=\beta_1(F)\wedge\cdots\wedge\beta_n(F)$,
where 
\[ \beta_i(F)=\left\{\begin{array}{cl}
\varepsilon_i(F)r,&\text{ if $i$ is an open root of $F$,}\\
r\omega_i,&\text{ if $i$ is a closed root of $F$,}\\
\sum_{\eta\in\mu_r}\eta^{\sum_{k\preceq i}ra_k}\omega_{i,j,\eta},
&\text{ if $i\to j$ in $F$.}
\end{array}\right. \]
The following can be proved by multiplying both sides by $z^r-w^r$:
\beq
\sum_{\eta\in\mu_r}\frac{\eta^a}{z-\eta w}=
\frac{rz^{a-1-r\lfloor\frac{a-1}{r}\rfloor}
w^{r-a+r\lfloor\frac{a-1}{r}\rfloor}}{z^r-w^r},\text{ for all $a\in\Z$.}
\eeq
Using this we obtain:
\bes
\begin{split}
\sum_{\eta\in\mu_r}\eta^a\omega_{i,j,\eta}
&=\left\{\begin{array}{cl}
\frac{rz_i^{a-r\lfloor\frac{a}{r}\rfloor}z_j^{r-a+r\lfloor\frac{a}{r}\rfloor}
(\omega_i-\omega_j)}{z_i^r-z_j^r},&\text{ if $r\nmid a$,}\\
\frac{r(z_i^{r-1}dz_i-z_j^{r-1}dz_j)}{z_i^r-z_j^r},&\text{ if $r\mid a$}
\end{array}\right.\\
&=\left\{\begin{array}{cl}
z_i^{a-r\lfloor\frac{a}{r}\rfloor}z_j^{-a+r\lceil\frac{a}{r}\rceil}
\varphi^*((z_i-z_j)^{-1}(\omega_i-\omega_j)),&\text{ if $r\nmid a$,}\\
\varphi^*(\omega_{i,j}),&\text{ if $r\mid a$.}
\end{array}\right.
\end{split}
\ees
So for any edge $i\to j$ in $F$,
\beq \label{secondbetaeqn}
\beta_i(F)=z_i^{\sum_{k\preceq i}ra_k-r\lfloor\sum_{k\preceq i}a_k\rfloor}
z_j^{-\sum_{k\preceq i}ra_k+r\lceil\sum_{k\preceq i}a_k\rceil}
\varphi^*(\bbeta_i(F)).
\eeq
If $i$ is a root of $F$, we have
\beq \label{firstbetaeqn}
\beta_i(F)=\varphi^*(\bbeta_i(F))=
z_i^{\sum_{k\preceq i}ra_k-r\lfloor\sum_{k\preceq i}a_k\rfloor}
\varphi^*(\bbeta_i(F)),
\eeq
since in this case $\sum_{k\preceq i}a_k\in\Z$ by assumption.
It only remains to note that for all $j$,
\bes
\sum_{k\preceq j}ra_k-r\lfloor\sum_{k\preceq j}a_k\rfloor-
\sum_{i\to j}\sum_{k\preceq i}ra_k+\sum_{i\to j}r\lceil
\sum_{k\preceq i}a_k\rceil
=ra_j+rb_j(F),
\ees
as required.
\epf

Combining Propositions \ref{cochainprop}, \ref{betaprop} and \ref{calcprop},
we deduce the following.
\bth \label{finalthm} Let $H^p(1,n;a_i)$ be the span
of the differential forms $\bbeta(F)$
for $F\in\bigcup_{k+l=p}\CalF(1,n;a_i)^{k,l}$.
\ben
\item $H^p(1,n;a_i)\subset\ker(d+\ell_{\sum a_i\omega_i}:\Omega^p(T(1,n))
\to\Omega^{p+1}(T(1,n)))$.
\item $H^p(1,n;a_i)\isomto H^p(T(1,n),\CalL_{\sum a_i\omega_i})$ 
via the map sending each
form to its class mod
$\im(d+\ell_{\sum a_i\omega_i}:\Omega^{p-1}(T(1,n))
\to\Omega^{p}(T(1,n)))$.
\item $H^p(1,n;a_i)$ has basis 
$\{\bbeta(F)\,|\,F\in\bigcup_{k+l=p}\CalF(1,n;a_i)^{\circ,k,l}\}$.
\item For $w\in Z_{(a_i)}$ and $F\in\CalF(1,n;a_i)^{k,l}$,
\[ w.\bbeta(F)=\varepsilon_n(w)\varepsilon(w,F)\bbeta(w.F). \]
\een
\eth
\bcr \label{finalcor}
We have
\[
\dim H^p(T(1,n),\CalL_{\sum a_i\omega_i})
=\sum_{k+l=p}|\CalF(1,n;a_i)^{\circ,k,l}|.
\]
In particular, $H^p(T(1,n),\CalL_{\sum a_i\omega_i})=0$ unless
$\sum_i a_i\in\Z$, and
\[ H^p(T(1,n),\CalL_{\frac{s}{n}\sum\omega_i})\neq 0\Longleftrightarrow
n-\gcd(s,n)\leq p\leq n. \]
\ecr
\bex
If all $a_i\in\Z$ and $r=1$, 
$\bbeta(F)=z_1^{-a_1}\cdots z_n^{-a_n}\alpha(F)$, and Theorem \ref{finalthm}
is the obvious translation 
of the results for $H^\bullet(T(1,n),\C)$ through the
isomorphism $\CalL_{\sum a_i\omega_i}\cong\C$.
\eex
\bex \label{mysteryex}
If all $a_i=\frac{s}{n}$ and $\gcd(s,n)=1$, 
$\CalF(1,n;\frac{s}{n})$ is the union of $\CalF(1,n;\frac{s}{n})^{n-1,0}$
and $\CalF(1,n;\frac{s}{n})^{n-1,1}$. The former set consists of trees
$T\in\CalT(1,n)$ with an open root; the second consists of the elements
$T^*$, which are the trees $T\in\CalT(1,n)$ with a closed root.
An easy calculation shows that for $T\in\CalT(1,n)^\circ$,
\[ \bbeta(T)=\frac{(-1)^{n-1}n\prod_{j=1}^n z_j^{b_j(T)-1}}{p_T}
\sum_i (-1)^{n-i}z_i \,dz_1\wedge\cdots\wedge
\widehat{dz_i}\wedge\cdots\wedge dz_n, \]
whereas
\[ \bbeta(T^*)=\frac{\prod_{j=1}^n z_j^{b_j(T)-1}}{p_T}
dz_1\wedge\cdots\wedge dz_n. \]
If $s=-1$ we have $b_j(T)=1$ for all $j$, so the factors
$\prod_{j=1}^n z_j^{b_j(T)-1}$ disappear.
(Compare the basis of $A^{n-1}(M(n))$ given by \eqref{nbceqn}.)
Note that this example can be handled by more general methods, since
$\frac{s}{n}\sum\omega_i-s\omega_n$ is not resonant (see Theorem \ref{esvthm}).
Indeed, Kawahara's result \cite[Theorem 2.1]{kawahara} implies that
$H^p(T(1,n),\CalL_{\frac{s}{n}\sum\omega_i})=0$
for $p<n-1$,
and \cite[Theorem 2.3]{kawahara} provides another basis of the non-vanishing
cohomologies.
\eex
\bex
Take $n=4$, all $a_i=-\frac{1}{2}$, $p=4$. The nine forests in $\bigcup_{k+l=4}
\CalF(1,4;-\frac{1}{2})^{\circ,k,l}$ are 
listed in the first column of the following
table (asterisks indicate closed roots).
The corresponding differential forms $\bbeta(F)$, giving the basis
of $H^4(1,4;-\frac{1}{2})$, are in the second 
column.
\[
\begin{array}{|c|c|}
\hline
1\to 2\to 3\to 4^* & \rule[-0.6cm]{0cm}{1.2cm}{\displaystyle
\frac{z_2}{(z_1-z_2)(z_2-z_3)(z_3-z_4)}}\,dz_1\wedge dz_2\wedge dz_3\wedge dz_4
\\\hline
1\to 2\to 4^*\leftarrow 3 & \rule[-0.6cm]{0cm}{1.2cm}{\displaystyle
\frac{z_2}{(z_1-z_2)(z_2-z_4)(z_3-z_4)}}\,dz_1\wedge dz_2\wedge dz_3\wedge dz_4
\\\hline
1\to 3\to 4^*\leftarrow 2 & \rule[-0.6cm]{0cm}{1.2cm}{\displaystyle
\frac{z_3}{(z_1-z_3)(z_2-z_4)(z_3-z_4)}}\,dz_1\wedge dz_2\wedge dz_3\wedge dz_4
\\\hline
1\to 4^*\leftarrow 3\leftarrow 2 & \rule[-0.6cm]{0cm}{1.2cm}{\displaystyle
\frac{z_3}{(z_1-z_4)(z_2-z_3)(z_3-z_4)}}\,dz_1\wedge dz_2\wedge dz_3\wedge dz_4
\\\hline
\begin{array}{ccccc}
1&\to&3&\to&4^*\\
&&\uparrow&&\\
&&2&&
\end{array} & {\displaystyle
\frac{z_3}{(z_1-z_3)(z_2-z_3)(z_3-z_4)}}\,dz_1\wedge dz_2\wedge dz_3\wedge dz_4
\\\hline
\begin{array}{ccccc}
1&\to&4^*&\leftarrow&3\\
&&\uparrow&&\\
&&2&&
\end{array} & {\displaystyle
\frac{z_4}{(z_1-z_4)(z_2-z_4)(z_3-z_4)}}\,dz_1\wedge dz_2\wedge dz_3\wedge dz_4
\\\hline
1\to 2^*\qquad 3\to 4^* & \rule[-0.6cm]{0cm}{1.2cm}{\displaystyle
\frac{1}{(z_1-z_2)(z_3-z_4)}}\,dz_1\wedge dz_2\wedge dz_3\wedge dz_4\\\hline
1\to 3^*\qquad 2\to 4^* & \rule[-0.6cm]{0cm}{1.2cm}{\displaystyle
\frac{1}{(z_1-z_3)(z_2-z_4)}}\,dz_1\wedge dz_2\wedge dz_3\wedge dz_4\\\hline
1\to 4^*\qquad 2\to 3^* & \rule[-0.6cm]{0cm}{1.2cm}{\displaystyle
\frac{1}{(z_1-z_4)(z_2-z_3)}}\,dz_1\wedge dz_2\wedge dz_3\wedge dz_4\\
\hline
\end{array}
\]
\eex

The wedge product of differential forms endows $H^\bullet(T(1,n),\C)$
with a ring structure, and $H^\bullet(T(1,n),\CalL_\omega)$ with the structure
of $H^\bullet(T(1,n),\C)$-module for any $\omega\in A^1(T(1,n))$. A
consequence of Theorem \ref{finalthm} is:
\bcr \label{modulecor}
The $H^\bullet(T(1,n),\C)$-module $H^\bullet(T(1,n),\CalL_{\sum a_i\omega_i})$
is generated by the images of those $\bbeta(F)$'s where
$F\in\CalF(1,n;a_i)^\circ$ has no closed roots or breakable
edges. In particular, $H^\bullet(T(1,n),\CalL_{\frac{s}{n}\sum\omega_i})$
is generated by $H^{n-\gcd(s,n)}(T(1,n),\CalL_{\frac{s}{n}\sum\omega_i})$.
\ecr
\bpf
For general $F\in\CalF(1,n;a_i)^\circ$, let 
$\widetilde{F}\in\CalF(1,n;a_i)^\circ$ be the forest obtained by deleting
all breakable edges and opening all closed roots. It is easy to see that
$b_j(F)=b_j(\widetilde{F})$ for all $j$, so $\bbeta(F)$ is, up to scalar,
the wedge product of $\bbeta(\widetilde{F})$ with a collection of $1$-forms
$\omega_i$ and $\omega_{i,j}$. Since these $1$-forms are closed, the first
statement follows. In the special case when all $a_i=\frac{s}{n}$,
it is clear that $\widetilde{F}\in\CalF(1,n;a_i)^{\circ,n-\gcd(s,n),0}$,
whence the second statement.
\epf
It would be interesting to find more information about this module,
such as a projective resolution.

Our remaining goal is to describe $H^\bullet(T(1,n),\CalL_{\sum a_i\omega_i})$
as a graded representation of $Z_{(a_i)}$. It is more convenient to describe
the isomorphic representation $A^\bullet(T(r,n))_{(a_i)}$; all results can be
translated easily using Propositions \ref{cochainprop} and \ref{calcprop}.

We start with the analogue of Proposition \ref{relprop}.
\bpr \label{newrelprop}
The following hold in $A^\bullet(T(r,n))_{(a_i)}$.
\ben
\item $\beta(F_1)+\beta(F_2)=\beta(F_3)$, whenever 
$F_1,F_2,F_3\in\CalF(1,n;a_i)$ are as in (1) of Proposition \ref{relprop}.
\item $\beta(F)+\beta(F')=0$, whenever $F,F'\in\CalF(1,n;a_i)$ are as in
(2) of Proposition \ref{relprop}.
\item $\beta(F)+\beta(F')=\beta(F'')$ whenever $F,F',F''\in\CalF(1,n;a_i)$
are as in (3) of Proposition \ref{relprop}.
\item $\beta(F)+\beta(F')=0$ whenever $F,F'\in\CalF(1,n;a_i)$ are as in
(3) of Proposition \ref{relprop} and the corresponding $F''$ is not in
$\CalF(1,n;a_i)$.
\een
\epr
\bpf
In each case (1)--(3), 
the translates of the forests by a fixed $\unz\in\mu_r^n$
are still a triple/pair of the right form, so these all follow instantly from
Proposition \ref{relprop}. The same applies in case (4), except that
$\beta(F'')=0$ by (1) of Proposition \ref{betaprop}.
\epf
The analogue of Theorem \ref{relthm}
follows, with an entirely analogous proof.
For any $F\in\CalF(1,n;a_i)^{k,l}$, let $\langle F\rangle$ denote the image of
$\beta(F)$ in 
\[ A^{k,l}(T(r,n))_{(a_i)}=(A^{k+l}(T(r,n))_l)_{(a_i)}/
(A^{k+l}(T(r,n))_{l+1})_{(a_i)}. \]
\bth \label{newrelthm}
With notation as above:
\ben
\item $A^p(T(r,n))_{(a_i)}$ can be defined 
abstractly as the vector space spanned
by $\{\beta(F)\,|\,F\in\bigcup_{k+l=p}\CalF(1,n;a_i)^{k,l}\}$ subject to the 
relations in Proposition \ref{newrelprop}.
\item $A^{k,l}(T(r,n))_{(a_i)}$ can be defined abstractly
as the vector space spanned by
$\{\lng F\rng\,|\,F\in\CalF(1,n;a_i)^{k,l}\}$ subject to the 
following relations:
\begin{itemize}
\item $\lng F_1\rng+\lng F_2\rng=\lng F_3\rng$ whenever 
$F_1,F_2,F_3\in\CalF(1,n;a_i)^{k,l}$ are as in
(1) of Proposition \ref{relprop};
\item $\lng F\rng+\lng F'\rng=0$ whenever 
$F,F'\in\CalF(1,n;a_i)^{k,l}$ are as in 
(2) or (3) of Proposition \ref{relprop}.
\end{itemize}
\item $A^{k,l}(T(r,n))_{(a_i)}$ has basis 
$\{\lng F\rng\,|\,F\in\CalF(1,n;a_i)^{\circ,k,l}\}$.
\item For $w\in Z_{(a_i)}$ and $F\in\CalF(1,n;a_i)^{k,l}$,
\[ w.\lng F\rng=\varepsilon_n(w)\,\varepsilon(w,F)\,\lng w.F\rng. \]
\item  As a representation of $Z_{(a_i)}$,
\[ A^p(T(r,n))_{(a_i)}\cong\bigoplus_{k+l=p}A^{k,l}(T(r,n))_{(a_i)}. \]
\een
\eth

For the remaining results we restrict to the case where all
$a_i=\frac{s}{r}$, $\gcd(s,r)=1$, and $r\mid n$. We 
have an analogue of Corollary \ref{bigindcor}, proved in
the same way (incorporating the extra feature that all trees must have a number
of vertices divisible by $r$):
\bcr \label{otherbigindcor}
As a representation of $S_n$, $\varepsilon_n\otimes A^{k,l}(T(r,n))_\sr$
is isomorphic to the following direct sum:
\[ \bigoplus_{\substack{\lambda^1,\lambda^2\\
|\lambda^1|+|\lambda^2|=n/r\\\ell(\lambda^1)=n-k-l\\
\ell(\lambda^2)=l}}\negthickspace\negthickspace
\Ind_{\substack{
((S_{r\lambda_1^1}\times\cdots\times S_{r\lambda_{n-k-l}^1})\phantom{blah}\\
\ \rtimes(S_{m_1(\lambda^1)}\times S_{m_2(\lambda^1)}\times\cdots))\\
\times((S_{r\lambda_1^2}\times\cdots\times S_{r\lambda_{l}^2})
\phantom{blahbla}\\
\
\rtimes(S_{m_1(\lambda^2)}\times S_{m_2(\lambda^2)}\times\cdots))}}^{S_n}
\negthickspace\negthickspace\negthickspace
\left(\negthickspace\begin{array}{c}
\varepsilon\otimes
\CalV{(1,r\lambda_1^1)}\otimes\cdots\otimes \CalV{(1,r\lambda_{n-k-l}^1)}\\
\phantom{\varepsilon}\otimes 
\CalV{(1,r\lambda_1^2)}\otimes\cdots\otimes \CalV{(1,r\lambda_{l}^2)}
\end{array}\negthickspace\right), \]
where $S_{r\lambda_a^j}$ acts on the $\CalV{(1,r\lambda_a^j)}$ factor,
$S_{m_i(\lambda^j)}$ acts by permuting the $\CalV{(1,r\lambda_a^j)}$
factors where $\lambda_a^j=i$, and $\varepsilon$ denotes
the product of the sign characters of the $S_{m_i(\lambda^1)}$ components.
\ecr
\bex
If $r=n$, Corollary \ref{otherbigindcor} says that $\varepsilon_n\otimes
A^{n-1,0}(T(n,n))_{(\frac{s}{n})}$ and $\varepsilon_n\otimes
A^{n-1,1}(T(n,n))_{(\frac{s}{n})}$ are isomorphic to $\CalV{(1,n)}$.
\eex
\bex \label{finalex}
In the case $n=4$, $r=2$, $s=1$, Corollary \ref{otherbigindcor} gives 
the following isomorphisms of representations of $S_4$:
\bes
\begin{split}
A^{2,0}(T(2,4))_{(\frac{1}{2})}&\cong\varepsilon_4\otimes 
\Ind_{(S_2\times S_2)\rtimes S_2}^{S_4}(\varepsilon\otimes
\CalV{(1,2)}\otimes\CalV{(1,2)})\cong V^{(31)},\\
A^{3,0}(T(2,4))_{(\frac{1}{2})}&\cong\varepsilon_4\otimes 
\Ind_{S_4}^{S_4}(\CalV{(1,4)})\cong V^{(31)}\oplus V^{(21^2)},\\
A^{2,1}(T(2,4))_{(\frac{1}{2})}&\cong\varepsilon_4\otimes 
\Ind_{S_2\times S_2}^{S_4}(\CalV{(1,2)}\otimes\CalV{(1,2)})\cong
V^{(4)}\oplus V^{(31)}\oplus V^{(2^2)},\\
A^{3,1}(T(2,4))_{(\frac{1}{2})}&\cong\varepsilon_4\otimes 
\Ind_{S_4}^{S_4}(\CalV{(1,4)})\cong V^{(31)}\oplus V^{(21^2)},\\
A^{2,2}(T(2,4))_{(\frac{1}{2})}&\cong\varepsilon_4\otimes 
\Ind_{(S_2\times S_2)\rtimes S_2}^{S_4}
(\CalV{(1,2)}\otimes\CalV{(1,2)})\cong V^{(4)}\oplus V^{(2^2)}.
\end{split}
\ees
Here the irreducible representations of $S_4$ are denoted $V^\lambda$, where
$\lambda$ is a partition of $4$; the convention is the one where
$V^{(4)}\cong 1$,
$V^{(1^4)}\cong\varepsilon_4$.
\eex
We can substitute Theorem \ref{lsthm} into Corollary \ref{otherbigindcor}, 
to obtain:
\bcr \label{othermainindcor}
As a representation of $S_n$, $\varepsilon_n\otimes A^{k,l}(T(r,n))_\sr$
is isomorphic to the following direct sum:
\[ \bigoplus_{\substack{\lambda^1,\lambda^2\\|\lambda^1|+|\lambda^2|=n/r\\
\ell(\lambda^1)=n-k-l\\\ell(\lambda^2)=l}}
\Ind_{\substack{((\mu_{r\lambda_1^1}\times\cdots\times 
\mu_{r\lambda_{n-k-l}^1})
\rtimes(S_{m_1(\lambda^1)}\times S_{m_2(\lambda^1)}\times\cdots))\\
\!\times((\mu_{\lambda_1^2}\times\cdots
\times\mu_{\lambda_{l}^2})
\rtimes(S_{m_1(\lambda^2)}\times S_{m_2(\lambda^2)}\times\cdots))}}^{S_n}
(\varepsilon\psi),
\]
where $\psi$ is the character which takes the product of the 
$\mu_{r\lambda_a^j}$ components, and $\varepsilon$ is the product of the sign
characters of the $S_{m_i(\lambda^1)}$ components.
\ecr
\noindent
Summing over $k$ and $l$, we get the desired 
expression for
$A^\bullet(T(r,n))_\sr\cong H^\bullet(T(r,n),\C)_\sr
\cong H^\bullet(T(1,n),\CalL_{\frac{s}{r}\sum\omega_i})$ 
as a direct sum of inductions of one-dimensional characters. 

Furthermore,
there is an obvious relationship
between Corollaries \ref{bigindcor}
and \ref{otherbigindcor}, which gives the following 
explanation of \eqref{intromysteryeqn}.
\bcr \label{mysterycor}
As a representation of $S_n$,
\[ \varepsilon_n\otimes A^{k,l}(T(r,n))_\sr\cong
\Ind_{W(r,n/r)}^{S_n}(\rmdet_{n/r}\otimes A^{k-n+n/r,l}(T(r,n/r))). \]
\ecr
\bpf
For the left-hand side, use the expression of Corollary \ref{otherbigindcor},
with the right-hand expression of Theorem \ref{mysterythm} substituted for
each of the $\CalV{(1,r\lambda_a^j)}$ factors. 
For the right-hand side, use the expression given by Corollary
\ref{bigindcor}. Comparing, we get the result.
\epf
\bibliographystyle{siam}

\end{document}